# A MIXED SINGULAR/SWITCHING CONTROL PROBLEM FOR A DIVIDEND POLICY WITH REVERSIBLE TECHNOLOGY INVESTMENT[1]

By Vathana Ly Vath, Huyên Pham and Stéphane Villeneuve

*Université d'Evry, Université Paris 7 and Université de Toulouse 1*

We consider a mixed stochastic control problem that arises in Mathematical Finance literature with the study of interactions between dividend policy and investment. This problem combines features of both optimal switching and singular control. We prove that our mixed problem can be decoupled in two pure optimal stopping and singular control problems. Furthermore, we describe the form of the optimal strategy by means of viscosity solution techniques and smooth-fit properties on the corresponding system of variational inequalities. Our results are of a quasi-explicit nature. From a financial viewpoint, we characterize situations where a firm manager decides optimally to postpone dividend distribution in order to invest in a reversible growth opportunity corresponding to a modern technology. In this paper a reversible opportunity means that the firm may disinvest from the modern technology and return back to its old technology by receiving some gain compensation. The results of our analysis take qualitatively different forms depending on the parameters values.

**1. Introduction.** Stochastic optimization problems that involve both bounded variation control and/or optimal switching are becoming timely problems in the applied probability literature and, more particularly, in Mathematical Finance. On one hand, the study of singular stochastic control problems in corporate Finance originates with the research on optimal dividend policy for a firm whose cash reserve follows a diffusion model; see Jeanblanc and Shiryaev [11] and Choulli, Taksar and Zhou [3]. On the other hand, the combined singular/stopping control problems have emerged in target tracking models (see Davis and Zervos [6] and Karatzas, Ocone,

Received May 2006; revised June 2007.
[1]Supported by the Europlace Institute.

Wang and Zervos [12]) as well as in Mathematical Finance from firm investment theory. For instance, Guo and Pham [10] have studied the optimal time to activate production and to control it by buying or selling capital, while Zervos [18] has applied this type of mixed problem in the field of real options theory. Finally, the theory of investment under uncertainty for a firm that can operate a production activity in different modes has led to optimal switching problems which have received a lot of attention in recent years from the applied mathematics community; see Brekke and Oksendal [2], Duckworth and Zervos [8] and Ly Vath and Pham [13].

In this paper we consider a combined stochastic control problem that has emerged in a recent paper by Décamps and Villeneuve [5] with the study of the interactions between dividend policy and investment under uncertainty. These authors have studied the interaction between dividend policy and irreversible investment decision in a growth opportunity. Our aim is to extend this work by relaxing the irreversible feature of the growth opportunity. In other words, we shall consider a firm with a technology in place that has the opportunity to invest in a new technology that increases its profitability. The firm self-finances the opportunity cost on its cash reserve. Once installed, the manager can decide to return back to the old technology by receiving some cash compensation. The mathematical formulation of this problem leads to a combined singular control/switching control for a one dimensional diffusion process. The diffusion process may take two regimes old or new that are switched at stopping times decisions. Within a regime, the manager has to choose a dividend policy that maximizes the expected value of all payouts until bankruptcy or regime transition. The transition from one regime to another incurs a cost or a benefit. The problem is to find the optimal mixed strategy that maximizes the expected returns.

Our analysis is rich enough to address several important questions that have arisen recently in the real option literature (see the book of Dixit and Pyndick [7] for an overview of this literature). What is the effect of financing constraints on investment decision? When is it optimal to postpone dividends distribution in order to invest? Basically, two assumptions in the real option theory are that the investment decision is made independently of the financial structure of the investment firm and also that the cash process generated by the investment is independent of any managerial decision. In contrast, our model studies the investment under uncertainty with the following set of assumptions. The firm is cash constrained and must finance its investments on its cash benefits, and the cash process generated by the investment depends only on the managerial decision to pay or not pay dividends, to quit or not quit the project. Our major finding is to characterize the natural intuition that the manager will delay dividend payments if the investment is sufficiently valuable.

As usual in stochastic control theory, the problem developed in this paper leads via the dynamic programming principle to a Hamilton–Jacobi–Bellman equation which forms in this paper a system of coupled variational inequalities. Therefore, a classical approach based on a verification theorem fails since it is very difficult to guess the shape of both the value function and the optimal strategy. To circumvent this difficulty, we use a viscosity solution approach and a uniqueness result combined with smooth-fit properties for determining the solution to the HJB system. As a by product, we also determine the shape of switching regions. Our findings take qualitatively different forms depending on both the profit rates of each technology and the transition costs.

The paper is organized as follows. We formulate the combined stochastic control problem in Section 2. In Section 3 we characterize by means of viscosity solutions, the system of variational inequalities satisfied by the value function, and we also state some regularity properties. Section 4 is devoted to qualitative results concerning the switching regions and in Section 5 we give the quasi-explicit computation and description of the value function and the optimal strategies.

**2. Model formulation: a mixed switching/singular control problem.** We consider a firm whose activities generate cash process. The manager of the firm acts in the best interest of its shareholders and maximizes the expected present value of dividends up to bankruptcy when the cash reserve becomes negative. The firm has at any time the possibility to invest in a modern technology that increases the drift of the cash from $\mu_0$ to $\mu_1$ without affecting the volatility $\sigma$. This growth opportunity requires a fixed cost $g > 0$ self-financed by the cash reserve. Moreover, we consider a reversible investment opportunity for the firm: the manager can decide to return back to the old technology by receiving some fixed gain compensation $(1 - \lambda)g$, with $0 < \lambda < 1$.

The mathematical formulation of this mixed singular/switching control problem is as follows. Let $W$ be a Brownian motion on a filtered probability space $(\Omega, \mathcal{F}, \mathbb{F} = (\mathcal{F}_t)_{t \geq 0}, \mathbb{P})$ satisfying the usual conditions.
- A strategy decision for the firm is a singular/switching control $\alpha = (Z, (\tau_n)_{n \geq 1}) \in \mathcal{A}$, where $Z \in \mathcal{Z}$, the set of $\mathbb{F}$-adapted cadlag nondecreasing processes, $Z_{0^-} = 0$, $(\tau_n)_n$ is an increasing sequence of stopping times, $\tau_n \to \infty$. $Z$ represents the total amount of dividends paid until time $t$, $(\tau_n)$ the switching technology (regimes) time decisions. By convention, regime $i = 0$ represents the old technology and $i = 1$ the modern technology.
- Starting from an initial state $(x, i) \in \mathbb{R} \times \{0, 1\}$ for the cash-regime value, and given a control $\alpha \in \mathcal{A}$, the dynamics of the cash reserve process of a firm is governed by

(2.1) $$dX_t = \mu_{I_t} dt + \sigma\, dW_t - dZ_t - dK_t, \qquad X_{0^-} = x,$$

where

$$I_t = \sum_{n\geq 0}(i\mathbf{1}_{\tau_{2n}\leq t<\tau_{2n+1}} + (1-i)\mathbf{1}_{\tau_{2n+1}\leq t<\tau_{2n+2}}), \qquad I_{0^-} = i,$$

(2.2)

$$K_t = \sum_{n\geq 0}(g_{i,1-i}\mathbf{1}_{\tau_{2n+1}\leq t<\tau_{2n+2}} + g_{1-i,i}\mathbf{1}_{\tau_{2n+2}\leq t<\tau_{2n+3}}),$$

with

$$0 \leq \mu_0 < \mu_1, \qquad \sigma > 0,$$
$$g_{01} = g > 0, \qquad g_{10} = -(1-\lambda)g < 0, \qquad 0 < \lambda < 1.$$

(Here we used the convention $\tau_0 = 0$.) We denote by $(X^{x,i}, I^i)$ the solution to (2.1)–(2.2) (as usual, we omit the dependence in the control $\alpha$ when there is no ambiguity). The time of strict bankruptcy is defined as

$$T = T^{x,i,\alpha} = \inf\{t \geq 0 : X_t^{x,i} < 0\},$$

and we set by convention $X_t^{x,i} = X_T^{x,i}$ for $t \geq T$. Thus, for $t \in [T \wedge \tau_{2n}, T \wedge \tau_{2n+1})$, the cash reserve $X^{x,i}$ is in technology $i$ (its drift term is $\mu_i$), while for $t \in [T \wedge \tau_{2n+1}, T \wedge \tau_{2n+2})$, $X^{x,i}$ is in technology $1-i$ (its drift term is $\mu_{1-i}$). Moreover,

$$X_{T\wedge\tau_{2n+1}}^{x,i} = X_{(T\wedge\tau_{2n+1})^-}^{x,i} - g_{i,1-i} \qquad \text{on } \{\tau_{2n+1} < T\},$$
$$X_{T\wedge\tau_{2n+2}}^{x,i} = X_{(T\wedge\tau_{2n+2})^-}^{x,i} - g_{1-i,i} \qquad \text{on } \{\tau_{2n+2} < T\}.$$

The optimal firm value is

(2.3) $$v_i(x) = \sup_{\alpha\in\mathcal{A}} \mathbb{E}\left[\int_0^{T^-} e^{-\rho t}\,dZ_t\right], \qquad x \in \mathbb{R}, i = 0, 1.$$

Here, we used the notation $\int_0^{T^-} e^{-\rho t}\,dZ_t = \int_{[0,T)} e^{-\rho t}\,dZ_t$. Notice that $v_i$ is nonnegative, and $v_i(x) = 0$ for $x < 0$. Since $T = T^{x,i,\alpha}$ is obviously nondecreasing in $x$, the value functions $v_i$ are clearly nondecreasing.

**3. Dynamic programming and general properties on the value functions.** We first introduce some notation. We denote by $R^{x,i}$ the cash reserve in absence of dividends distribution and in regime $i$, that is, the solution to

(3.1) $$dR_t^{x,i} = \mu_i\,dt + \sigma\,dW_t, \qquad R_0^{x,i} = x.$$

The associated second-order differential operator is denoted $\mathcal{L}_i$:

$$\mathcal{L}_i\varphi(x) = \mu_i\varphi'(x) + \tfrac{1}{2}\sigma^2\varphi''(x).$$

In view of the dynamic programming principle, recalled below [see (3.20)], we formally expect that the value functions $v_i$, $i = 0, 1$, satisfy the system of variational inequalities:

$$\min[\rho v_i(x) - \mathcal{L}_i v_i(x), v_i'(x) - 1, v_i(x) - v_{1-i}(x - g_{i,1-i})] = 0, \quad (3.2)$$
$$x > 0, i = 0, 1.$$

This statement will later be proved rigorously by means of viscosity solutions. For the moment, we first state a standard comparison principle for this system of PDE.

PROPOSITION 3.1. *Suppose that $\varphi_i$, $i = 0, 1$, are two smooth functions on $(0, \infty)$ s.t. $\varphi_i(0^+) := \lim_{x \downarrow 0} \varphi_i(x) \geq 0$, and*

$$\min[\rho \varphi_i(x) - \mathcal{L}_i \varphi_i(x), \varphi_i'(x) - 1, \varphi_i(x) - \varphi_{1-i}(x - g_{i,1-i})] \geq 0, \quad (3.3)$$
$$x > 0, \ i = 0, 1,$$

*where we set by convention $\varphi_i(x) = 0$ for $x < 0$. Then, we have $v_i \leq \varphi_i$, $i = 0, 1$.*

PROOF. Given an initial state-regime value $(x, i) \in (0, \infty) \times \{0, 1\}$, take an arbitrary control $\alpha = (Z, (\tau_n), n \geq 1) \in \mathcal{A}$, and set for $m > 0$, $\theta_{m,n} = \inf\{t \geq T \wedge \tau_{2n} : X_t^{x,i} \geq m \text{ or } X_t^{x,i} \leq 1/m\} \nearrow \infty$ a.s. when $m$ goes to infinity. Apply then Itô's formula to $e^{-\rho t} \varphi_i(X_t^{x,i})$ between the stopping times $T \wedge \tau_{2n}$ and $\tau_{m,2n+1} := T \wedge \tau_{2n+1} \wedge \theta_{m,n}$. Notice that for $T \wedge \tau_{2n} \leq t < \tau_{m,2n+1}$, $X_t^{x,i}$ stays in regime $i$. Then, we have

$$e^{-\rho \tau_{m,2n+1}} \varphi_i(X_{\tau_{m,2n+1}^-}^{x,i})$$
$$= e^{-\rho(T \wedge \tau_{2n})} \varphi_i(X_{T \wedge \tau_{2n}}^{x,i}) + \int_{T \wedge \tau_{2n}}^{\tau_{m,2n+1}} e^{-\rho t}(-\rho \varphi_i + \mathcal{L}_i \varphi_i)(X_t^{x,i}) \, dt$$
$$(3.4)$$
$$+ \int_{T \wedge \tau_{2n}}^{\tau_{m,2n+1}} e^{-\rho t} \sigma \varphi_i'(X_t^{x,i}) \, dW_t - \int_{T \wedge \tau_{2n}}^{\tau_{m,2n+1}} e^{-\rho t} \varphi_i'(X_t^{x,i}) \, dZ_t^c$$
$$+ \sum_{T \wedge \tau_{2n} \leq t < \tau_{m,2n+1}} e^{-\rho t} [\varphi_i(X_t^{x,i}) - \varphi_i(X_{t^-}^{x,i})],$$

where $Z^c$ is the continuous part of $Z$. We make the convention that when $T \leq \tau_n$, $(T \wedge \theta)^- = T$ for all stopping time $\theta > \tau_n$ a.s., so that (3.4) holds true a.s. for all $n, m$ [recall that $\varphi_i(X_T^{x,i}) = 0$]. Since $\varphi_i' \geq 1$, we have by the mean-value theorem $\varphi_i(X_t^{x,i}) - \varphi_i(X_{t^-}^{x,i}) \leq X_t^{x,i} - X_{t^-}^{x,i} = -(Z_t - Z_{t^-})$ for $T \wedge \tau_{2n} \leq t < \tau_{m,2n+1}$. By using also the supersolution inequality of $\varphi_i$, taking expectation in the above Itô's formula, and noting that the integrand in the stochastic integral term is bounded by a constant (depending on $m$), we have

$$\mathbb{E}[e^{-\rho \tau_{m,2n+1}} \varphi_i(X^{x,i}_{\tau^-_{m,2n+1}})]$$
$$\leq \mathbb{E}[e^{-\rho(T \wedge \tau_{2n})} \varphi_i(X^{x,i}_{T \wedge \tau_{2n}})] - \mathbb{E}\left[\int_{T \wedge \tau_{2n}}^{\tau_{m,2n+1}} e^{-\rho t} \, dZ^c_t\right]$$
$$- \mathbb{E}\left[\sum_{T \wedge \tau_{2n} \leq t < \tau_{m,2n+1}} e^{-\rho t}(Z_t - Z_{t^-})\right],$$

and so

$$\mathbb{E}[e^{-\rho(T \wedge \tau_{2n})} \varphi_i(X^{x,i}_{T \wedge \tau_{2n}})] \geq \mathbb{E}\left[\int_{T \wedge \tau_{2n}}^{\tau^-_{m,2n+1}} e^{-\rho t} \, dZ_t + e^{-\rho \tau_{m,2n+1}} \varphi_i(X^{x,i}_{\tau^-_{m,2n+1}})\right].$$

By sending $m$ to infinity, with Fatou's lemma, we obtain

(3.5)
$$\mathbb{E}[e^{-\rho(T \wedge \tau_{2n})} \varphi_i(X^{x,i}_{T \wedge \tau_{2n}})]$$
$$\geq \mathbb{E}\left[\int_{T \wedge \tau_{2n}}^{(T \wedge \tau_{2n+1})^-} e^{-\rho t} \, dZ_t + e^{-\rho(T \wedge \tau_{2n+1})} \varphi_i(X^{x,i}_{(T \wedge \tau_{2n+1})^-})\right].$$

Now, as $\varphi_i(x) \geq \varphi_{1-i}(x - g_{i,1-i})$ and recalling $X^{x,i}_{T \wedge \tau_{2n+1}} = X^{x,i}_{(T \wedge \tau_{2n+1})^-} - g_{i,1-i}$ on $\{\tau_{2n+1} < T\}$, we have

(3.6)
$$\varphi_i(X^{x,i}_{(T \wedge \tau_{2n+1})^-}) \geq \varphi_{1-i}(X^{x,i}_{(T \wedge \tau_{2n+1})^-} - g_{i,1-i})$$
$$= \varphi_{1-i}(X^{x,i}_{(T \wedge \tau_{2n+1})}) \qquad \text{on } \{\tau_{2n+1} < T\}.$$

Moreover, notice that $\varphi_i$ is nonnegative as $\varphi_i(0^+) \geq 0$ and $\varphi'_i \geq 1$. Hence, since $\varphi_{1-i}(X^{x,i}_{(T \wedge \tau_{2n+1})}) = \varphi_{i-1}(X^{x,i}_T) = 0$ on $\{T \leq \tau_{2n+1}\}$, we see that inequality (3.6) also holds on $\{T \leq \tau_{2n+1}\}$ and so a.s., therefore, plugging into (3.5), we have

$$\mathbb{E}[e^{-\rho(T \wedge \tau_{2n})} \varphi_i(X^{x,i}_{T \wedge \tau_{2n}})]$$
$$\geq \mathbb{E}\left[\int_{T \wedge \tau_{2n}}^{(T \wedge \tau_{2n+1})^-} e^{-\rho t} \, dZ_t + e^{-\rho(T \wedge \tau_{2n+1})} \varphi_{1-i}(X^{x,i}_{T \wedge \tau_{2n+1}})\right].$$

Similarly, we have from the supersolution inequality of $\varphi_{1-i}$

$$\mathbb{E}[e^{-\rho(T \wedge \tau_{2n+1})} \varphi_{1-i}(X^{x,i}_{T \wedge \tau_{2n+1}})]$$
$$\geq \mathbb{E}\left[\int_{T \wedge \tau_{2n+1}}^{(T \wedge \tau_{2n+2})^-} e^{-\rho t} \, dZ_t + e^{-\rho(T \wedge \tau_{2n+2})} \varphi_i(X^{x,i}_{T \wedge \tau_{2n+2}})\right].$$

By iterating these two previous inequalities for all $n$, we then obtain

$$\varphi_i(x) \geq \mathbb{E}\bigg[\int_0^{(T\wedge\tau_{2n})^-} e^{-\rho t}\, dZ_t + e^{-\rho(T\wedge\tau_{2n})}\varphi_i(X^{x,i}_{T\wedge\tau_{2n}})\bigg],$$

$$\geq \mathbb{E}\bigg[\int_0^{(T\wedge\tau_{2n})^-} e^{-\rho t}\, dZ_t\bigg] \qquad \forall n \geq 0,$$

since $\varphi_i$ is nonnegative. By sending $n$ to infinity, we obtain the required result from the arbitrariness of the control $\alpha$. $\square$

As a corollary, we show a linear growth condition on the value functions.

COROLLARY 3.1. *We have*

(3.7) $\quad v_0(x) \leq x + \dfrac{\mu_1}{\rho}, \qquad v_1(x) \leq x + \dfrac{\mu_1}{\rho} + (1-\lambda)g, \qquad x > 0.$

PROOF. We set $\varphi_0(x) = x + \frac{\mu_1}{\rho}$, $\varphi_1(x) = x + \frac{\mu_1}{\rho} + (1-\lambda)g$, on $(0,\infty)$, and $\varphi_i(x) = 0$ for $x < 0$. A straightforward computation shows that we have the supersolution properties for $\varphi_i$, $i = 0,1$:

$$\min[\rho\varphi_0(x) - \mathcal{L}_0\varphi_0(x), \varphi_0'(x) - 1, \varphi_0(x) - \varphi_1(x-g)] \geq 0, \qquad x > 0,$$

$$\min[\rho\varphi_1(x) - \mathcal{L}_1\varphi_1(x), \varphi_1'(x) - 1, \varphi_1(x) - \varphi_0(x + (1-\lambda)g)] \geq 0, \qquad x > 0.$$

We then conclude from Proposition 3.1. $\square$

The next result states the initial-boundary data for the value functions.

PROPOSITION 3.2. (1) *The value function $v_0$ is continuous on $(0,\infty)$ and satisfies*

(3.8) $\qquad\qquad v_0(0^+) := \lim_{x\downarrow 0} v_0(x) = 0.$

(2) *The value function $v_1$ satisfies*

(3.9) $\qquad\qquad v_1(0^+) := \lim_{x\downarrow 0} v_1(x) = v_0((1-\lambda)g).$

PROOF. (1) (a) We first state (3.8). For $x > 0$, let us consider the drifted Brownian $R^{x,1}$, defined in (3.1), and denote $\theta_0 = \inf\{t \geq 0 \colon R^{x,1}_t = 0\}$. It is well known that

(3.10) $\qquad\qquad \mathbb{E}\bigg[\sup_{0\leq t\leq \theta_0} R^{x,1}_t\bigg] \to 0 \qquad \text{as } x \downarrow 0.$

We also have

(3.11) $\qquad\qquad \sup_{0\leq t\leq \theta_0} R^{x,1}_t \downarrow 0 \qquad \text{a.s. as } x \downarrow 0.$

Fix some $r > 0$, and denote $\theta_r = \inf\{t \geq 0 : R_t^{x,1} = r\}$. It is also well known that

$$\mathbb{P}[\theta_0 > \theta_r] \to 0 \qquad \text{as } x \downarrow 0. \tag{3.12}$$

Let $\alpha = (Z, (\tau_n)_{n \geq 1})$ be an arbitrary policy in $\mathcal{A}$, and denote $\eta = T \wedge \theta_r = T^{x,0,\alpha} \wedge \theta_r$. Since $\mu_0 < \mu_1$ and $g_{01} > 0$, $g_{01} + g_{10} > 0$, we notice that $X_t^{x,0} \leq R_t^{x,1} - Z_t \leq R_t^{x,1}$ for all $t \geq 0$. Hence, $T \leq \theta_0$, $Z_t \leq R_t^{x,1}$ for $t < T$, and, in particular, $Z_{\eta^-} \leq R_\eta^{x,1}$. We then write

$$
\begin{aligned}
\mathbb{E}\left[\int_0^{T^-} e^{-\rho t} \, dZ_t\right] &= \mathbb{E}\left[\int_0^{\eta^-} e^{-\rho t} \, dZ_t\right] + \mathbb{E}\left[1_{T > \eta} \int_\eta^{T^-} e^{-\rho t} \, dZ_t\right] \\
&\leq \mathbb{E}[Z_{\eta^-}] + \mathbb{E}\left[\mathbb{E}\left[1_{T > \eta} \int_\eta^{T^-} e^{-\rho t} \, dZ_t \Big| \mathcal{F}_{\theta_r^-}\right]\right] \\
&\leq \mathbb{E}[R_\eta^{x,1}] + \mathbb{E}\left[1_{T > \theta_r} \mathbb{E}\left[\int_{\theta_r}^{T^-} e^{-\rho t} \, dZ_t \Big| \mathcal{F}_{\theta_r^-}\right]\right] \\
&\leq \mathbb{E}[R_\eta^{x,1}] + \mathbb{E}[1_{T > \theta_r} e^{-\rho \theta_r} v_0(X_{\theta_r^-}^{x,0})],
\end{aligned}
\tag{3.13}
$$

where we also used in the second inequality the fact that on $\{T > \eta\}$, $\eta = \theta_r$, and $\theta_r$ is a predictable stopping time, and in the last inequality the definition of the value function $v_0$. Now, since $v_0$ is nondecreasing, we have $v_0(X_{\theta_r^-}^{x,0}) \leq v_0(r)$. Moreover, recalling that $T \leq \theta_0$, inequality (3.13) yields

$$0 \leq v_0(x) \leq \mathbb{E}\left[\sup_{0 \leq t \leq \theta_0} R_t^{x,1}\right] + v_0(r) \mathbb{P}[\theta_0 > \theta_r] \longrightarrow 0 \qquad \text{as } x \downarrow 0, \tag{3.14}$$

from (3.10)–(3.12). This proves $v_0(0^+) = 0$.

(b) We next prove the continuity of $v_0$ at any $y > 0$. Let $\alpha = (Z, (\tau_n)_{n \geq 1}) \in \mathcal{A}$, $X^{y,0}$ be the corresponding process and $T = T^{y,0,\alpha}$ its bankruptcy time. According to (3.10) and (3.12), given a fixed $r > 0$, for any arbitrary small $\varepsilon > 0$, one can find $0 < \delta < y$ s.t. for $0 < x < \delta$,

$$\mathbb{E}\left[\sup_{0 \leq t \leq \theta_0} R_t^{x,1}\right] + v_0(r) \mathbb{P}[\theta_0 > \theta_r] \leq \varepsilon.$$

Then, following the same lines of proof as for (3.13)–(3.14), we show

$$\mathbb{E}\left[\int_\theta^{T^-} e^{-\rho t} \, dZ_t\right] \leq \varepsilon, \tag{3.15}$$

for any $0 < x < \delta$ and stopping time $\theta$ s.t. $X_\theta^{y,0} \leq x$. Given $0 < x < \delta$, consider the state process $X^{y-x,0}$ starting from $y - x$ in regime 0, and controlled by $\alpha$. Denote $\theta$ its bankruptcy time, that is, $\theta = T^{y-x,0,\alpha} = \inf\{t \geq 0 : X_t^{y-x,0} < 0\}$. Notice that $X_t^{y-x,0} = X_t^{y,0} - x$ for $t \leq \theta \leq T$, and so

$$X_\theta^{y,0} = X_\theta^{y-x,0} + x \leq x.$$

From (3.15), we then have

$$\mathbb{E}\left[\int_0^{T^-} e^{-\rho t}\, dZ_t\right] = \mathbb{E}\left[\int_0^{\theta^-} e^{-\rho t}\, dZ_t\right] + \mathbb{E}\left[\int_\theta^{T^-} e^{-\rho t}\, dZ_t\right]$$
$$\leq v_0(y-x) + \varepsilon.$$

From the arbitrariness of $\alpha$, and recalling that $v_0$ is nondecreasing, this implies

$$0 \leq v_0(y) - v_0(y-x) \leq \varepsilon,$$

which shows the continuity of $v_0$.

(2) Given an arbitrary control $\alpha = (Z, (\tau_n)_{n\geq 1}) \in \mathcal{A}$, let us consider the control $\tilde\alpha = (\tilde Z, (\tilde\tau_n)_{n\geq 1}) \in \mathcal{A}$ defined by $\tilde Z = Z$, $\tilde\tau_1 = 0$, $\tilde\tau_n = \tau_{n-1}$, $n \geq 2$. Then, for all $x > 0$, and by stressing the dependence of the state process on the control, we have $X_t^{x,1,\tilde\alpha} = X_t^{x+(1-\lambda)g,0,\alpha}$ for $0 \leq t < T^{x,1,\tilde\alpha} = T^{x+(1-\lambda)g,0,\alpha}$. We deduce

$$v_1(x) \geq \mathbb{E}\left[\int_0^{(T^{x,1,\tilde\alpha})^-} e^{-\rho t}\, d\tilde Z_t\right] = \mathbb{E}\left[\int_0^{(T^{x+(1-\lambda)g,0,\alpha})^-} e^{-\rho t}\, dZ_t\right],$$

which implies, from the arbitrariness of $\alpha$,

(3.16) $\qquad v_1(x) \geq v_0(x + (1-\lambda)g), \qquad x > 0.$

On the other hand, starting in the regime $i = 1$, for $x \geq 0$, let $\alpha = (Z, (\tau_n)_{n\geq 1})$ be an arbitrary control in $\mathcal{A}$. We denote $T_1 = T \wedge \tau_1 = T^{x,1,\alpha} \wedge \tau_1$, and we write

(3.17) $\quad \mathbb{E}\left[\int_0^{T^-} e^{-\rho t}\, dZ_t\right] = \mathbb{E}\left[\int_0^{T_1^-} e^{-\rho t}\, dZ_t\right] + \mathbb{E}\left[\mathbf{1}_{T>\tau_1}\int_{\tau_1}^{T^-} e^{-\rho t}\, dZ_t\right].$

The first term in the r.h.s. of (3.17) is dealt similarly as in (3.13)–(3.14): we set $\eta_1 = T_1 \wedge \theta_r$ with $\theta_r = \inf\{t \geq 0 : R_t^{x,1} = r\}$ for some fixed $r > 0$, and we notice that $X_t^{x,1} = R_t^{x,1} - Z_t \leq R_t^{x,1}$ for $t < \tau_1$. Hence, $T_1 \leq \theta_0 = \inf\{t \geq 0 : R_t^{x,1} = 0\}$, and $Z_{\eta_1^-} \leq R_{\eta_1}^{x,1} \leq \sup_{0 \leq t \leq \theta_0} R_t^{x,1}$. Then, as in (3.13)–(3.14), we have

(3.18) $\qquad \mathbb{E}\left[\int_0^{T_1^-} e^{-\rho t}\, dZ_t\right] \leq \mathbb{E}\left[\sup_{0\leq t\leq \theta_0} R_t^{x,1}\right] + v_1(r)\mathbb{P}[\theta_0 > \theta_r].$

For the second term in the r.h.s. of (3.17), since there is a change of regime at $\tau_1$ from $i = 1$ to $i = 0$, and by definition of the value function $v_0$, we have

$$\mathbb{E}\left[\mathbf{1}_{T>\tau_1}\int_{\tau_1}^{T^-} e^{-\rho t}\, dZ_t\right] = \mathbb{E}\left[\mathbf{1}_{T>\tau_1}\mathbb{E}\left[\int_{\tau_1}^{T^-} e^{-\rho t}\, dZ_t\Big|\mathcal{F}_{\tau_1}\right]\right]$$
$$\leq \mathbb{E}[\mathbf{1}_{T>\tau_1} e^{-\rho\tau_1} v_0(X_{\tau_1}^{x,1})]$$

(3.19)
$$\leq \mathbb{E}[1_{T > \tau_1} v_0(X^{x,1}_{\tau_1^-} + (1-\lambda)g)]$$
$$\leq \mathbb{E}\left[v_0\left(\sup_{0 \leq t \leq \theta_0} R^{x,1}_t + (1-\lambda)g\right)\right].$$

Here, we used in the second inequality the fact that $X^{x,1}_{\tau_1} = X^{x,1}_{\tau_1^-} + (1-\lambda)g$ on $\{\tau_1 < T\}$, and in the last one the observation that $X^{x,1}_t \leq R^{x,1}_t$ for $t < \tau_1$, and $\tau_1 = T_1 \leq \theta_0$ on $\{\tau_1 < T\}$. Hence, by combining (3.16)–(3.19), we obtain

$$v_0(x + (1-\lambda)g) \leq v_1(x) \leq \mathbb{E}\left[\sup_{0 \leq t \leq \theta_0} R^{x,1}_t\right] + v_1(r)\mathbb{P}[\theta_0 > \theta_r]$$
$$+ \mathbb{E}\left[v_0\left(\sup_{0 \leq t \leq \theta_0} R^{x,1}_t + (1-\lambda)g\right)\right].$$

Finally, by using the continuity of $v_0$, the limits (3.10)–(3.12), as well as the linear growth condition (3.7) of $v_0$, which allows to apply the dominated convergence theorem, we conclude that $v_1(0^+) = v_0((1-\lambda)g)$. □

REMARK 3.1. There is some asymmetry between the two value functions $v_0$ and $v_1$. Actually, $v_0$ is continuous at 0: $v_0(0^+) = v_0(0^-) = 0$, while it is not the case for $v_1$, since $v_1(0^+) = v_0((1-\lambda)g) > 0 = v_1(0^-)$. When the reserve process in regime 0 approaches zero, we are ineluctably absorbed by this threshold. On the contrary, in regime 1, when the reserve process approaches zero, we have the possibility to change the regime, which pushes us above the bankruptcy threshold by receiving $(1-\lambda)g$. In particular, at this stage, we do not know yet the continuity of $v_1$ on $(0, \infty)$. This will be proved in Theorem 3.1 as a consequence of the dynamic programming principle. In the sequel, we set by convention $v_i(0) = v_i(0^+)$ for $i = 0, 1$.

We shall assume that the following dynamic programming principle holds: for any $(x, i) \in \mathbb{R}_+ \times \{0, 1\}$, we have

(DP) $$v_i(x) = \sup_{\alpha \in \mathcal{A}} \mathbb{E}\left[\int_0^{(T \wedge \theta \wedge \tau_1)^-} e^{-\rho t} dZ_t\right.$$

(3.20) $$+ e^{-\rho(T \wedge \theta \wedge \tau_1)}(v_i(X^{x,i}_{T \wedge \theta}) 1_{T \wedge \theta < \tau_1}$$
$$\left. + v_{1-i}(X^{x,i}_{\tau_1}) 1_{\tau_1 \leq T \wedge \theta})\right],$$

where $\theta$ is any stopping time, possibly depending on $\alpha \in \mathcal{A}$ in (3.20).

We then have the PDE characterization of the value functions $v_i$.

THEOREM 3.1. *The value functions $v_i$, $i = 0, 1$, are continuous on $(0, \infty)$, and are the unique viscosity solutions with linear growth condition on $(0, \infty)$ and boundary data $v_0(0) = 0$, $v_1(0) = v_0((1 - \lambda)g)$ to the system of variational inequalities:*

$$\min[\rho v_i(x) - \mathcal{L}_i v_i(x), v_i'(x) - 1, v_i(x) - v_{1-i}(x - g_{i,1-i})] = 0,$$
(3.21)
$$x > 0, i = 0, 1.$$

Actually, we prove some more regularity results on the value functions.

PROPOSITION 3.3. *The value functions $v_i$, $i = 0, 1$, are $C^1$ on $(0, \infty)$. Moreover, if we set, for $i = 0, 1$,*

$$\mathcal{S}_i = \{x \geq 0 : v_i(x) = v_{1-i}(x - g_{i,1-i})\},$$
$$\mathcal{D}_i = \{x > 0 : v_i'(x) = 1\},$$
$$\mathcal{C}_i = (0, \infty) \setminus (\mathcal{S}_i \cup \mathcal{D}_i),$$

*then $v_i$ is $C^2$ on the open set $\mathcal{C}_i \cup \operatorname{int}(\mathcal{D}_i)$ of $(0, \infty)$, and we have in the classical sense*

$$\rho v_i(x) - \mathcal{L}_i v_i(x) = 0, \qquad x \in \mathcal{C}_i.$$

REMARK 3.2. From the variational inequality (3.21), and since the value functions $v_i$, $i = 0, 1$, are $C^1$ on $(0, \infty)$, we have $v_i' \geq 1$, which implies, in particular, that $v_i$ is strictly increasing on $(0, \infty)$.

The proofs of Theorem 3.1 and Proposition 3.3 follow and combine essentially arguments from [10] for singular control, and [14] for switching control, and are postponed to Appendix A and B.

$\mathcal{S}_i$ is the switching region from technology $i$ to $1 - i$, $\mathcal{D}_i$ is the dividend region in technology $i$, and $\mathcal{C}_i$ is the continuation region in technology $i$. Notice from the boundary conditions on $v_i$ that $\mathcal{S}_i$ may contain 0. We denote $\mathcal{S}_i^* = \mathcal{S}_i \setminus \{0\}$.

## 4. Qualitative results on the switching regions.

4.1. *Benchmarks.* We consider the firm value without investment/disinvestment in technology $i = 0$:

(4.1) $$\hat{V}_0(x) = \sup_{Z \in \mathcal{Z}} \mathbb{E}\left[\int_0^{T_0^-} e^{-\rho t} \, dZ_t\right],$$

where $T_0 = \inf\{t \geq 0 : X_t \leq 0\}$ is the time bankruptcy of the cash reserve in regime 0:

$$dX_t = \mu_0 \, dt + \sigma \, dW_t - dZ_t, \qquad X_{0^-} = x.$$

By convention, we set $\hat{V}_0(x) = 0$ for $x < 0$. It is known that $\hat{V}_0$, as the value function of a pure singular control problem, is characterized as the unique continuous viscosity solution on $(0, \infty)$, with linear growth condition to the variational inequality

$$\text{(4.2)} \qquad \min[\rho \hat{V}_0 - \mathcal{L}_0 \hat{V}_0, \hat{V}_0' - 1] = 0, \qquad x > 0,$$

and boundary data

$$\hat{V}_0(0) = 0.$$

Actually, $\hat{V}_0$ is $C^2$ on $(0, \infty)$ and explicit computations of this standard singular control problem are developed in Shreve, Lehoczky and Gaver [16], Jeanblanc and Shiryaev [11], or Radner and Shepp [15]:

$$\hat{V}_0(x) = \begin{cases} \dfrac{f_0(x)}{f_0'(\hat{x}_0)}, & 0 \leq x \leq \hat{x}_0, \\ x - \hat{x}_0 + \dfrac{\mu_0}{\rho}, & x \geq \hat{x}_0, \end{cases}$$

where

$$f_0(x) = e^{m_0^+ x} - e^{m_0^- x}, \qquad \hat{x}_0 = \frac{1}{m_0^+ - m_0^-} \ln\left(\frac{(m_0^+)^2}{(m_0^-)^2}\right),$$

and $m_0^- < 0 < m_0^+$ are roots of the characteristic equation

$$\rho - \mu_0 m - \tfrac{1}{2}\sigma^2 m^2 = 0.$$

In other words, this means that the optimal cash reserve process is given by the reflected diffusion process at the threshold $\hat{x}_0$ with an optimal dividend process given by the local time at this boundary. When the firm starts with a cash reserve $x \geq \hat{x}_0$, the optimal dividend policy is to distribute immediately the amount $x - \hat{x}_0$ and then follow the dividend policy characterized by the local time.

As a second benchmark, we consider the firm value problem in technology $i = 1$ with nonnegative constant liquidation value $L$ to be fixed later:

$$w_1^L(x) = \sup_{Z \in \mathcal{Z}} \mathbb{E}\left[\int_0^{T_1^-} e^{-\rho t}\, dZ_t + e^{-\rho T_1} L\right],$$

$T_1 = \inf\{t \geq 0 : X_t \leq 0\}$ is the time bankruptcy of the cash reserve in regime 1:

$$dX_t = \mu_1\, dt + \sigma\, dW_t - dZ_t, \qquad X_{0^-} = x.$$

By convention, we set $w_1^L(x) = 0$ for $x < 0$. Again, as value function of a pure singular control problem, $w_1^L$ is characterized as the unique continuous viscosity solution on $(0, \infty)$, with linear growth condition to the variational inequality

$$(4.3) \qquad \min[\rho w_1^L - \mathcal{L}_1 w_1^L, (w_1^L)' - 1] = 0, \qquad x > 0,$$

and boundary data

$$(4.4) \qquad w_1^L(0) = L.$$

Actually, $w_1^L$ is $C^2$ on $(0, \infty)$ and explicit computations of this singular control problem are developed in Boguslavskaya [1]:
- If $L \geq \frac{\mu_1}{\rho}$, then

$$w_1^L(x) = x + L, \qquad x \geq 0.$$

The optimal strategy is to distribute the initial cash reserve immediately, and so to liquidate the firm at $X_t = 0$ by changing of technology to regime $i = 0$ and receiving $L$.
- If $L < \frac{\mu_1}{\rho}$, then

$$(4.5) \qquad w_1^L(x) = \begin{cases} \dfrac{1 - L h_1'(\hat{x}_1)}{f_1'(\hat{x}_1)} f_1(x) + L h_1(x), & 0 \leq x \leq x_1^L, \\ x - x_1^L + \dfrac{\mu_1}{\rho}, & x \geq x_1^L, \end{cases}$$

with

$$f_1(x) = e^{m_1^+ x} - e^{m_1^- x}, \qquad h_1(x) = e^{m_1^- x},$$

$m_1^- < 0 < m_1^+$, the roots of the characteristic equation

$$\rho - \mu_1 m - \tfrac{1}{2}\sigma^2 m^2 = 0,$$

and $x_1^L$ the solution to

$$(4.6) \qquad L \frac{h_1(x) f_1'(x) - h_1'(x) f_1(x)}{f_1'(x)} + \frac{f_1(x)}{f_1'(x)} = \frac{\mu_1}{\rho}.$$

The optimal cash reserve process is given by the reflected diffusion process at the threshold $x_1^L$ with an optimal dividend process given by the local time at this boundary. When the firm starts with a cash reserve $x \geq x_1^L$, the optimal dividend policy is to distribute immediately the amount $x - x_1^L$ and then follow the dividend policy characterized by the local time. In the sequel we shall denote

$$\hat{V}_1 = w_1^L \quad \text{and} \quad \hat{x}_1 = x_1^L \qquad \text{when } L = \hat{V}_0((1-\lambda)g).$$

$L = \hat{V}_0((1-\lambda)g)$ is the minimal received liquidation value when one switches to regime 0 at $x = 0$ and does not switch anymore.

REMARK 4.1. It is known (see, e.g., [1]) that $\hat{V}_0$ and $w_1^L$ are concave on $(0, \infty)$. As a consequence, $\hat{V}_0$ and $w_1^L$ are globally Lipschitz since their first derivatives are bounded near zero.

REMARK 4.2. We have $v_0 \geq \hat{V}_0$ and $v_1 \geq \hat{V}_1$ on $(0, \infty)$. This is rather clear since the class of controls over which maximization is taken in $\hat{V}_0$ and $\hat{V}_1$ is included in the class of controls of $v_0$ and $v_1$. This may be justified more rigorously by a maximum principle argument and by noting that $v_0$ and $v_1$ are (viscosity) supersolution to the variational inequality satisfied respectively by $\hat{v}_0$ and $\hat{V}_1$, with the same boundary data.

We first show the intuitive result that the value function for the dividend policy problem is nondecreasing in the rate of return of the cash reserve.

LEMMA 4.1.
$$\hat{V}_1(x) \geq \hat{V}_0(x + (1-\lambda)g) \qquad \forall x \geq 0.$$

PROOF. We set $w_1(x) = \hat{V}_1(x - (1-\lambda)g)$ for $x \geq (1-\lambda)g$. From (4.3), we see that $\hat{w}_1$ satisfies on $[(1-\lambda)g, \infty)$
$$w_1'(x) = \hat{V}_1'(x - (1-\lambda)g) \geq 1,$$
$$(\rho w_1 - \mathcal{L}_0 w_1)(x) = (\rho - \mathcal{L}_1 \hat{V}_1 + (\mu_1 - \mu_0)\hat{V}_1')(x - (1-\lambda)g) > 0,$$
since $\mu_1 > \mu_0$ and $\hat{V}_1$ is increasing. Moreover, $w_1((1-\lambda)g) = \hat{V}_1(0) = \hat{V}_0((1-\lambda)g)$. By the standard maximum principle on the variational inequality (4.2), we deduce that $w_1 \geq \hat{V}_0$ on $[(1-\lambda)g, \infty)$, which implies the required result. □

The next result precises conditions under which the value function in the old technology is larger than the value function in the modern technology after paying the switching cost from the old to the modern regimes.

LEMMA 4.2. *Suppose that* $\hat{V}_0((1-\lambda)g) < \frac{\mu_1}{\rho}$. *Then,*
$$\hat{V}_0(x) \geq \hat{V}_1(x - g) \qquad \forall x \geq 0 \quad \text{if and only if} \quad \frac{\mu_1 - \mu_0}{\rho} \leq \hat{x}_1 + g - \hat{x}_0.$$

PROOF. Similar arguments as in Lemma 2.1 in Decamps and Villeneuve [5]. □

REMARK 4.3. Recalling that $\hat{V}_0$ and $\hat{V}_1$ are increasing and concave, the above lemma shows also that if $\frac{\mu_1 - \mu_0}{\rho} > \hat{x}_1 + g - \hat{x}_0$, then there exists $\hat{x}_{01} \geq g$ s.t.

$$\max(\hat{V}_0(x), \hat{V}_1(x-g)) = \begin{cases} \hat{V}_0(x), & x \leq \hat{x}_{01}, \\ \hat{V}_1(x-g), & x > \hat{x}_{01}. \end{cases}$$

4.2. *Preliminary results on the switching regions.* In this section we shall state some preliminary qualitative results concerning the switching regions.

LEMMA 4.3. *If $x \in \mathcal{S}_i$, then $x - g_{i,1-i} \notin \mathcal{S}_{1-i}$.*

PROOF. Since $v_i(x) > v_i(x - \lambda g)$ for every $x > 0$ and $i \in \{0,1\}$, we have, for $x \in \mathcal{S}_i$,

$$v_{1-i}(x - g_{i,1-i}) = v_i(x) > v_i(x - \lambda g) = v_i(x - g_{i,1-i} - g_{1-i,i}).$$

Therefore, $x - g_{i,1-i} \notin \mathcal{S}_{1-i}$ for $x \in \mathcal{S}_i$. □

Let us recall the notation $\mathcal{S}_i^* = \mathcal{S}_i \setminus \{0\}$. We have the following inclusion:

LEMMA 4.4. $\mathcal{S}_1^* \subset \mathcal{D}_1$.

PROOF. We make a proof by contradiction by assuming that there exists some $x \in \mathcal{S}_1^* \setminus \mathcal{D}_1$. According to Proposition 3.3, we have $v_0'(x + (1-\lambda)g) = v_1'(x) > 1$, and so $x + (1-\lambda)g \notin \mathcal{D}_0$. Applying Lemma 4.3 with $i = 1$ implies $x + (1-\lambda)g \in \mathcal{C}_0$. Therefore,

$$\begin{aligned}
\rho v_1(x) - \mathcal{L}_1 v_1(x) &= \rho v_1(x) - \mathcal{L}_0 v_1(x) + (\mu_0 - \mu_1) v_1'(x) \\
&= \rho v_0(x + (1-\lambda)g) - \mathcal{L}_0 v_0(x + (1-\lambda)g) \\
&\quad + (\mu_0 - \mu_1) v_1'(x) \\
&= (\mu_0 - \mu_1) v_1'(x) \quad \text{since } x + (1-\lambda)g \in \mathcal{C}_0 \\
&< 0,
\end{aligned}$$

which contradicts Theorem 3.1. □

We now introduce the following definition.

DEFINITION 4.1. *$y$ is a left boundary of the closed set $\mathcal{D}_i$ if there is some $\delta > 0$ such that $y - \varepsilon$ does not belong to $\mathcal{D}_i$ for every $0 < \varepsilon < \delta$.*

LEMMA 4.5.  *Let $y > 0$ be a left boundary of $\mathcal{D}_i$:*

- *If there is some $\varepsilon > 0$ such that $(y - \varepsilon, y) \subset \mathcal{C}_i$, then $v_i(y) = \frac{\mu_i}{\rho}$.*
- *If not, $v_i(y) = \frac{\mu_{1-i}}{\rho}$.*

PROOF. Since $y$ is a left boundary of $\mathcal{D}_i$, there is some $\varepsilon > 0$ such that $(y - \varepsilon, y) \subset \mathcal{C}_i \cup \mathcal{S}_i$. Therefore, two cases have to be considered:

*Case* 1: If $(y - \varepsilon, y) \subset \mathcal{C}_i$, then, according to Proposition 3.3, $v_i$ is twice differentiable at $x$, for $y - \varepsilon < x < y$, and satisfies $v_i'(x) = 1$ and $v_i''(x) = 0$. Therefore, we have

$$0 = \rho v_i(x) - \mathcal{L}_i v_i(x) = \rho v_i(x) - \mu_i v_i'(x) - \frac{\sigma^2}{2} v_i''(x).$$

By sending $x$ to $y$, we obtain that $v_i(y) = \frac{\mu_i}{\rho}$.

*Case* 2: If not, there is an increasing sequence $(y_n)_n$ valued in $\mathcal{S}_i$, and converging to $y$ which therefore belongs to $\mathcal{S}_i$. We then have $v_i(y_n) = v_{1-i}(y_n - g_{i,1-i})$ and also $v_i'(y_n) > 1$ for $n$ great enough since $y$ is a left boundary of $\mathcal{D}_i$. Thus, $y_n - g_{i,1-i} \notin \mathcal{D}_{1-i}$. Moreover, according to Lemma 4.3, we also have $y_n - g_{i,1-i} \notin \mathcal{S}_{1-i}$ and, therefore, $y_n - g_{i,1-i} \in \mathcal{C}_{1-i}$ or, equivalently,

$$\rho v_{1-i}(y_n - g_{i,1-i}) - \mathcal{L}_{1-i} v_{1-i}(y_n - g_{i,1-i}) = 0.$$

By letting $n$ tends to $\infty$, we obtain $v_{1-i}(y - g_{i,1-i}) = \frac{\mu_{1-i}}{\rho}$. Since $y \in \mathcal{S}_i$, this implies $v_i(y) = v_{1-i}(y - g_{i,1-i}) = \frac{\mu_{1-i}}{\rho}$. □

The next result shows that the switching region from modern technology $i = 1$ to the old technology $i = 0$ is either reduced to the zero threshold or to the entire state reserve domain $\mathbb{R}_+$, depending on the gain $(1 - \lambda)g$ for switching from regime 1 to regime 0.

PROPOSITION 4.1. *The two following cases arise:*

  (i) *If $v_0((1 - \lambda)g) < \frac{\mu_1}{\rho}$, then $\mathcal{S}_1 = \{0\}$.*
  (ii) *If $v_0((1 - \lambda)g) \geq \frac{\mu_1}{\rho}$, then $\mathcal{S}_1 = \mathcal{D}_1 = \mathbb{R}_+$.*

PROOF. (i) Assume $v_0((1 - \lambda)g) < \frac{\mu_1}{\rho}$. We shall make a proof by contradiction by considering the existence of some $x_0 \in \mathcal{S}_1^*$. By Lemma 4.4, one can introduce the finite nonnegative number

$$\underline{x} = \inf\{y > 0 : [y, x_0] \subset \mathcal{D}_1\}.$$

Hence, $\underline{x}$ is a left boundary of $\mathcal{D}_1$. Moreover, Lemma 4.5 gives $v_1(\underline{x}) = \frac{\mu_1}{\rho}$ or $\frac{\mu_0}{\rho}$.

1. We first check that $\underline{x} > 0$. If not, we would have $v_1(y) = y + v_0((1-\lambda)g)$ for any $0 < y < x_0$. But, in this case, we have, for $0 < y < x_0$,
$$\rho v_1(y) - \mathcal{L}_1 v_1(y) = \rho(y + v_0((1-\lambda)g)) - \mu_1.$$
Therefore, under the assumption (i), $\rho v_1(y) - \mathcal{L}_1 v_1(y) < 0$ for $y$ small enough which is a contradiction.

2. We now prove that $v_1(\underline{x}) = \frac{\mu_1}{\rho}$. To see this, we shall show that the closed set $\mathcal{D}_1$ is an interval of $\mathbb{R}_+$. Letting $a, b \in \mathcal{D}_1$ with $a < b$, we want to show that $(a, b) \subset \mathcal{D}_1$. If not, from Lemma 4.4, we can find a subinterval $(c, d)$ with $c, d \in \mathcal{D}_1$ and $(c, d) \subset \mathcal{C}_1$. But, for $c < x < d$, we have
$$0 = \rho v_1(x) - \mathcal{L}_1 v_1(x) = \rho v_1(x) - \mu_1 v_1'(x) - \frac{\sigma^2}{2} v_1''(x).$$
By sending $x$ to $c$ and $d$, we obtain that $v_1(c) = v_1(d) = \frac{\mu_1}{\rho}$, which contradicts the fact that $v_1$ is strictly increasing. Since $\mathcal{D}_1$ is an interval of $\mathbb{R}_+$, we have $\underline{x} = \inf \mathcal{D}_1$. Thus, recalling that $\underline{x} > 0$, we can find, from Lemma 4.4, some $\varepsilon > 0$ such that $(\underline{x} - \varepsilon, \underline{x}) \subset \mathcal{C}_1$, and deduce from Lemma 4.5 that $v_1(\underline{x}) = \frac{\mu_1}{\rho}$.

3. We now introduce
$$\bar{x} = \inf\{y \geq \underline{x} | y \in \mathcal{S}_1\}.$$
Observe that $\bar{x} + (1-\lambda)g \in \mathcal{D}_0$. Moreover, according to Lemma 4.3, $\bar{x} + (1-\lambda)g \notin \mathcal{S}_0$ and, thus, a left neighborhood of $\bar{x} + (1-\lambda)g$ belongs to $\mathcal{C}_0$. We first prove that $\bar{x} + (1-\lambda)g$ cannot be a left boundary of $\mathcal{D}_0$. On the contrary, we would have, from Lemma 4.5,
$$v_1(\bar{x}) = v_0(\bar{x} + (1-\lambda)g) = \frac{\mu_0}{\rho} < \frac{\mu_1}{\rho} = v_1(\underline{x}),$$
which contradicts the fact that $v_1$ is increasing. Therefore, $\bar{x} + (1-\lambda)g \in \overset{o}{\mathcal{D}}_0$, and we can find $y < \bar{x}$ such that $y + (1-\lambda)g$ is a left boundary of $\mathcal{D}_0$. Hence,
$$v_1(\bar{x}) = v_0(\bar{x} + (1-\lambda)g) = \bar{x} - y + v_0(y + (1-\lambda)g) \leq \bar{x} - y + v_1(y).$$
Since the reverse inequality is always true, we obtain that $y \in \mathcal{S}_1$, which contradicts the definition of $\bar{x}$. We conclude that $\bar{x}$ cannot be strictly positive, which is a contradiction with the first step. This proves finally that $x_0 \in \mathcal{S}_1^*$.

(ii) Assume that $v_0((1-\lambda)g) \geq \frac{\mu_1}{\rho}$. Let $y$ be a left boundary of $\mathcal{D}_1$. We shall prove that $y$ necessarily equals zero. If not, according to Lemma 4.5, $v_1(y) \leq \frac{\mu_1}{\rho} \leq v_1(0)$, where the second inequality comes from the hypothesis and (3.9). Since the function $v_1$ is strictly increasing, we get the desired contradiction. Therefore, $\mathcal{D}_1 = [0, a]$. It remains to prove that $a$ is infinite.

From Lemma 4.4, the open set $(a, \infty)$ belongs to $\mathcal{C}_1$ if $a < \infty$. Using the regularity of $v_1$ on $\mathcal{C}_1$, we get by the same reasoning as in the proof of Lemma 4.5 that $v_1(a) = \frac{\mu_1}{\rho}$, which gives the same contradiction as before. Hence, $\mathcal{D}_1 = [0, \infty)$. We then have, for any $x > 0$,

$$v_1(x) = x + v_0((1-\lambda)g) \leq v_0(x + (1-\lambda)g).$$

Since the reverse inequality is always true by definition, we conclude that $\mathcal{S}_1 = [0, \infty)$. □

The next proposition describes the structure of the switching region from technology $i = 0$ to $i = 1$, in the case where the growth rate $\mu_1$, in the modern technology $i = 1$, is large enough.

PROPOSITION 4.2. *Suppose that*

$$\frac{\mu_1 - \mu_0}{\rho} > \hat{x}_1 + g - \hat{x}_0 \quad and \quad \hat{V}_0((1-\lambda)g) < \frac{\mu_1}{\rho}.$$

*Then, there exists $x_{01}^* \in [g, \infty)$ s.t.*

$$\mathcal{S}_0^* = [x_{01}^*, \infty).$$

PROOF. We first notice that $\mathcal{S}_0^* \neq \varnothing$. On the contrary, we would have $v_0 = \hat{V}_0$, and so $\hat{V}_0(x) \geq v_1(x - g) \geq \hat{V}_1(x - g)$ for all $x$, which is in contradiction with Lemma 4.2. Moreover, since $v_1(x - g) = v_0(x) > 0$ for all $x \in \mathcal{S}_0^*$, we deduce that $\mathcal{S}_0^* \subset [g, \infty)$ and so

$$x_{01}^* := \inf \mathcal{S}_0^* \in [g, \infty).$$

Let us now consider the function

$$w_0(x) = \begin{cases} v_0(x), & x < x_{01}^*, \\ v_1(x - g), & x \geq x_{01}^*. \end{cases}$$

We claim that $w_0$ is a viscosity solution, with linear growth condition and boundary data $w_0(0^+) = 0$, to

$$\min[\rho w_0(x) - \mathcal{L}_0 w_0(x), w_0'(x) - 1, w_0(x) - v_1(x - g)] = 0, \qquad x > 0.$$

For $x < x_{01}^*$, this is clear since $w_0 = v_0$ on $(0, x_{01}^*)$. For $x > x_{01}^*$, we see that $w_0' \geq 1$ and

$$\rho w_0 - \mathcal{L}_0 w_0 = (\rho v_1 - \mathcal{L}_1 v_1 + (\mu_1 - \mu_0) v_1')(x - g)$$
$$\geq (\mu_1 - \mu_0) v_1'(x - g) \geq 0.$$

Hence, the viscosity property is also satisfied for $x > x_{01}^*$. It remains to check the viscosity property for $x = x_{01}^*$. The viscosity subsolution property at $x_{01}^*$ is trivial since $w_0(x_{01}^*) = v_1(x_{01}^* - g)$. For the viscosity supersolution property, take some $C^2$ test function $\varphi$ s.t. $x_{01}^*$ is a local minimum of $w_0 - \varphi$. From the smooth-fit condition of the value function $v_0$ at the switching boundary, it follows that $w_0$ is $C^1$ at $x_{01}^*$. Hence, $w_0'(x_{01}^*) = \varphi'(x_{01}^*)$. Moreover, since $w_0 = v_0$ is $C^2$ for $x < x_{01}^*$, we also have $\varphi''(x_{01}^*) \leq w_0''(x_{01}^{*-}) := \lim_{x \nearrow x_{01}^*} w''(x)$. Since $\rho w_0(x) - \mathcal{L}_0 w_0(x) \geq 0$ for $x < x_{01}^*$, we deduce by sending $x$ to $x_{01}^*$:

$$\rho w_0(x_{01}^*) - \mathcal{L}_0 \varphi(x_{01}^*) \geq 0.$$

This implies the required viscosity supersolution inequality at $x = x_{01}^*$. By uniqueness, we conclude that $w_0 = v_0$, which proves that $\mathcal{S}_0^* = [x_{01}^*, \infty)$. $\square$

**5. Main result and description of the solution.** We give an explicit description of the structure of the solution to our control problem, which depends crucially on parameter values.

5.1. *The case:* $\hat{V}_0((1-\lambda)g) \geq \frac{\mu_1}{\rho}$.

THEOREM 5.1. *Suppose that* $\hat{V}_0((1-\lambda)g) \geq \frac{\mu_1}{\rho}$. *Then, we have* $v_0(x) = \hat{V}_0(x)$ *and* $v_1(x) = \hat{V}_0(x + (1-\lambda)g) = x + (1-\lambda)g - x_0 + \frac{\mu_0}{\rho}$. *It is optimal to never switch from regime* 0 *to regime* 1. *In regime* 1, *it is optimal to distribute all the surplus as dividends and to switch to regime* 0.

PROOF. Under the condition of the theorem, and since $v_0 \geq \hat{V}_0$, we have $v_0((1-\lambda)g) \geq \frac{\mu_1}{\rho}$. By Proposition 4.1, this implies $\mathcal{S}_1 = \mathcal{D}_1 = \mathbb{R}_+$. Recalling also the boundary data $v_1(0) = v_0((1-\lambda)g)$, we get $v_1(x) = x + v_0((1-\lambda)g)$ for $x \geq 0$. We next prove that the region $\mathcal{S}_0^*$ is empty. To see this, we have to prove that, for $x \geq g$, $v_0(x) \geq v_1(x-g)$. Let us consider for $x \geq g$ the function $\theta(x) = v_0(x) - (x - g + v_0((1-\lambda)g))$. Since $\lambda > 0$, we have $\theta(g) > 0$. Moreover, $\theta'(x) = v_0'(x) - 1 \geq 0$. Thus, $\theta(x) > 0$ for $x \geq g$, which is equivalent to $\mathcal{S}_0^* = \varnothing$. As a consequence, $v_0$ is a smooth solution of the variational inequality

$$\min[\rho v(x) - \mathcal{L}_0 v(x), v'(x) - 1] = 0,$$

with initial condition $v(0) = 0$. By uniqueness, we deduce that $v_0 = \hat{V}_0$. To close the proof, it suffices to note that $\hat{V}_0((1-\lambda)g) \geq \frac{\mu_1}{\rho}$ implies that $(1-\lambda)g \geq \hat{x}_0$. Therefore, $v_0((1-\lambda)g) = (1-\lambda)g - x_0 + \frac{\mu_0}{\rho}$. $\square$

5.2. *The case:* $\hat{V}_0((1-\lambda)g) < \frac{\mu_1}{\rho}$. First observe that, in this case, we have

$$v_0((1-\lambda)g) < \frac{\mu_1}{\rho}.$$

Indeed, on the contrary, from Theorem 5.1, we would get $v_0 = \hat{V}_0$, and so an obvious contradiction $\hat{V}_0((1-\lambda)g) \geq \frac{\mu_1}{\rho}$ with the considered case. From Proposition 4.1, we then have $\mathcal{S}_1 = \{0\}$ so that $v_1$ is the unique viscosity solution to

$$\min[\rho v_1 - \mathcal{L}_1 v_1, v_1' - 1] = 0, \qquad x > 0,$$

with the boundary data $v_1(0) = v_0((1-\lambda)g)$. Therefore, $v_1$ is the firm value problem in technology $i = 1$ with liquidation value $v_0((1-\lambda)g)$:

$$(5.1) \qquad v_1(x) = \sup_{Z \in \mathcal{Z}} \mathbb{E}\left[\int_0^{T_1^-} e^{-\rho t}\, dZ_t + e^{-\rho T_1} v_0((1-\lambda)g)\right].$$

The form of $v_1$ is described in (4.5) with liquidation value $L = v_0((1-\lambda)g)$: we denote $x_1 = x_1^L$ the corresponding threshold.

REMARK 5.1. Since $v_1$ and $\hat{V}_1$ are increasing with $v_1(x_1) = \hat{V}_1(\hat{x}_1) = \frac{\mu_1}{\rho}$, we have $x_1 \leq \hat{x}_1$.

Notice that the expression of $v_1$ is not completely explicit since we do not know at this stage the liquidation value $v_0((1-\lambda)g)$. The next result give an explicit solution when

$$\frac{\mu_1 - \mu_0}{\rho} \leq \hat{x}_1 + g - \hat{x}_0.$$

THEOREM 5.2. *Suppose that*

$$(5.2) \qquad \hat{V}_0((1-\lambda)g) < \frac{\mu_1}{\rho} \leq \frac{\mu_0}{\rho} + \hat{x}_1 + g - \hat{x}_0.$$

*Then $v_0 = \hat{V}_0$ and $v_1 = \hat{V}_1$. It is never optimal, once in regime $i = 0$, to switch to regime $i = 1$. In regime 1, it is optimal to switch to regime 0 at the threshold $x = 0$.*

PROOF. From Lemma 4.1 and Lemma 4.2, and recalling the variational inequalities (4.2) and (4.3), we see that $\hat{V}_0$ and $\hat{V}_1$ are viscosity solutions to

$$\min[\rho \hat{V}_0(x) - \mathcal{L}_0 \hat{V}_0(x), \hat{V}_0'(x) - 1, \hat{V}_0(x) - \hat{V}_1(x-g)] = 0, \qquad x > 0,$$

$$\min[\rho \hat{V}_1(x) - \mathcal{L}_1 \hat{V}_1(x), \hat{V}_1'(x) - 1, \hat{V}_1(x) - \hat{V}_0(x+(1-\lambda)g)] = 0, \qquad x > 0,$$

together with the boundary data $V_0(0^+) = 0$ and $\hat{V}_1(0^+) = \hat{V}_0((1-\lambda)g)$. By uniqueness to this system of variational inequalities, we conclude that $(v_0, v_1) = (\hat{V}_0, \hat{V}_1)$. □

In the sequel, we suppose that

(5.3) $$\frac{\mu_1 - \mu_0}{\rho} > \hat{x}_1 + g - \hat{x}_0.$$

From Proposition 4.2, the switching region from regime 0 to regime 1 has the form

$$\mathcal{S}_0^* = \{x > 0 : v_0(x) = v_1(x - g)\} = [x_{01}^*, \infty),$$

for some $x_{01}^* \in [g, \infty)$. Moreover, since $x_1 \leq \hat{x}_1$ (see Remark 5.1), the above condition (5.3) implies $\frac{\mu_1 - \mu_0}{\rho} > x_1 + g - \hat{x}_0$. By the same arguments as in Remark 4.3, there exists some $\bar{x}_{01} \geq g$ s.t.

$$\max(\hat{V}_0(x), v_1(x-g)) = \begin{cases} \hat{V}_0(x), & x \leq \bar{x}_{01}, \\ v_1(x-g), & x > \bar{x}_{01}. \end{cases}$$

Following [5], we introduce the pure stopping time problem

(5.4) $$\bar{v}_0(x) = \sup_{\tau \in \mathcal{T}} \mathbb{E}[e^{-\rho(\tau \wedge T_0)} \max(\hat{V}_0(R_{\tau \wedge T_0}^{x,0}), v_1(R_{\tau \wedge T_0}^{x,0} - g))],$$

where $\mathcal{T}$ denotes the set of stopping times valued in $[0, \infty]$. We also denote $\mathcal{E}_0$ the exercise region for $\bar{v}_0$:

$$\mathcal{E}_0 = \{x \geq 0 : \bar{v}_0(x) = \max(\hat{V}_0(x), v_1(x-g))\}.$$

The next result shows that the original mixed singular/switching control problems may be reformulated as a coupled pure optimal stopping time and pure singular problem.

THEOREM 5.3. *Suppose that*

(5.5) $$\hat{V}_0((1-\lambda)g) < \frac{\mu_1}{\rho} \quad \text{and} \quad \frac{\mu_1 - \mu_0}{\rho} > \hat{x}_1 + g - \hat{x}_0.$$

*Then, we have*

$$v_0 = \bar{v}_0$$

*and $v_1$ given by (5.1). Moreover,*

$$\mathcal{E}_0 = \{0 \leq x < \bar{x}_{01} : v_0(x) = \hat{V}_0(x)\} \cup [x_{01}^*, \infty).$$

PROOF. The proof follows along the lines of those of Theorem 3.1 in [5]. We will give only the road map of it in our context and omit the details.

Let us first note that the process $(e^{-\rho(t \wedge T_0)} v_0(R_{t \wedge T_0}^{x,0}))_{t \geq 0}$ is a supermartingale that dominates the function $\max(\hat{V}_0, v_1(\cdot - g))$. Therefore, according to the Snell envelope theory, we have $v_0 \geq \bar{v}_0$.

To prove the reverse inequality, it is enough to show that $\bar{v}_0' \geq 1$ (see Proposition 3.4 in [5]) and to use the uniqueness result of Theorem 3.1. To this end, we will precise the shape of the exercise region $\mathcal{E}_0$. According to Lemma 4.3 by Villeneuve [17], $\bar{x}_{01}$ does not belong to $\mathcal{E}_0$. Thus, the exercise region can be decomposed into two subregions

$$\mathcal{E}_{00} = \{x < \bar{x}_{01} : v_0(x) = \hat{V}_0(x)\}$$

and

$$\mathcal{E}_{01} = \{x > \bar{x}_{01} : v_0(x) = v_1(x - g)\}.$$

Two cases have to be considered:

*Case* (i). If the subregion $\mathcal{E}_{00}$ is empty, the optimal stopping problem defined by $\bar{v}_0$ can be solved explicitly, and we have (see [5], Lemma 3.3)

$$\bar{v}_0 = \begin{cases} \dfrac{e^{m_0^+ x} - e^{m_0^- x}}{e^{m_0^+ x_{01}^*} - e^{m_0^- x_{01}^*}} v_1(x_{01}^* - g), & x < x_{01}^*, \\ v_1(x - g), & x \geq x_{01}^*. \end{cases}$$

The smooth-fit principle allows us to conclude that $\bar{v}_0' \geq 1$ since $v_1' \geq 1$.

*Case* (ii). If the subregion $\mathcal{E}_{00}$ is nonempty, we can prove using the arguments of Proposition 3.5 and Lemma 3.4 in [5] that

$$\mathcal{E}_0 = [0, a] \cup [x_{01}^*, \infty),$$

with $a \geq \hat{x}_0$ and the value function $\bar{v}_0$ satisfies

$$\bar{v}_0(x) = A e^{m_0^+ x} + B e^{m_0^- x} \qquad \text{for } x \in (a, x_{01}^*).$$

The smooth-fit principle gives $\bar{v}_0'(a) = \hat{V}_0'(a) \geq 1$ and $\bar{v}_0'(x_{01}^*) = v_1'(x_{01}^* - g) \geq 1$. Clearly, $\bar{v}_0$ is convex in a right neighborhood of $a$ since $\hat{V}_0$ is linear at $a$. Therefore, if $\bar{v}_0$ remains convex on $(a, x_{01}^*)$, the proof is over. If not, the second derivative of $\bar{v}_0$ given by $A(m_0^+)^2 e^{m_0^+ x} + B(m_0^-)^2 e^{m_0^- x}$ vanishes at most one time on $(a, x_{01}^*)$, say, in $d$. Hence,

$$1 = \bar{v}_0'(a) \leq (\bar{v}_0)'(x) \leq \bar{v}_0'(d) \qquad \text{for } x \in (a, d)$$

and

$$1 \leq \bar{v}_0'(x_{01}^*) \leq \bar{v}_0'(x) \leq \bar{v}_0'(d) \qquad \text{for } x \in (d, x_{01}^*),$$

which completes the proof. □

Notice that the representation (5.1)–(5.4) of pure optimal singular and stopping problems for $v_1$ and $v_0$ is coupled, and so not easily computable. We decouple this representation by considering the sequence of pure optimal stopping and singular control problems, starting from $\hat{V}_1^{(0)} = \hat{V}_1$ and $\hat{V}_0^{(0)} = \hat{V}_0$:

$$\hat{V}_0^{(k)}(x) = \sup_{\tau \in \mathcal{T}} \mathbb{E}[e^{-\rho(\tau \wedge T_0)} \max(\hat{V}_0(R^{x,0}_{\tau \wedge T_0}), \hat{V}_1^{(k-1)}(R^{x,0}_{\tau \wedge T_0} - g))], \qquad k \geq 1,$$

$$\hat{V}_1^{(k)}(x) = \sup_{Z \in \mathcal{Z}} \mathbb{E}\left[\int_0^{T_1^-} e^{-\rho t}\, dZ_t + e^{-\rho T_1} \hat{V}_0^{(k)}((1-\lambda)g)\right], \qquad k \geq 1.$$

The next result shows the convergence of this procedure.

PROPOSITION 5.1. *Under the conditions* (5.5) *of Theorem* 5.3, *we have, for all* $x > 0$,

$$\lim_{k \to \infty} \hat{V}_0^{(k)}(x) = v_0(x), \qquad \lim_{k \to \infty} \hat{V}_1^{(k)}(x) = v_1(x).$$

PROOF. We will first prove that the increasing sequence $(\hat{V}_0^{(k)}, \hat{V}_1^{(k)})$ converges uniformly on every compact subset of $\mathbb{R}_+$. To see this, we will apply the Arzela–Ascoli theorem by first proving the equi-continuity of the functions $\hat{V}_i^{(k)}$. Let us first remark that the functions $\hat{V}_1^{(k)}$ are Lipschitz continuous uniformly in $k$ since they are concave with bounded first derivative (see Remark 4.2) independently of $k$. Let us also check that the functions $\hat{V}_0^{(k)}$ are Lipschitz continuous uniformly in $k$. Using the inequality $\max(a,b) - \max(c,d) \leq \max(a-c, b-d)$, and by setting

$$\Delta(x,y) = \max(\hat{V}_0(R^{x,0}_{\tau \wedge T_0}) - \hat{V}_0(R^{y,0}_{\tau \wedge T_0}), \hat{V}_1^{(k-1)}(R^{x,0}_{\tau \wedge T_0} - g) - \hat{V}_1^{(k-1)}(R^{y,0}_{\tau \wedge T_0} - g)),$$

we get by recalling also that $\hat{V}_0$ is Lipschitz (see Remark 4.2)

$$\begin{aligned}|\hat{V}_0^{(k)}(x) - \hat{V}_1^{(k)}(y)| &\leq \sup_{\tau \in \mathcal{T}} \mathbb{E}[e^{-\rho(\tau \wedge T_0)} |\Delta(x,y)|] \\ &\leq K_0 \sup_{\tau \in \mathcal{T}} \mathbb{E}[e^{-\rho(\tau \wedge T_0)} |R^{x,0}_{\tau \wedge T_0} - R^{y,0}_{\tau \wedge T_0}|] \\ &\leq K_0 |x-y| \sup_{\tau \in \mathcal{T}} \mathbb{E}[e^{-\rho(\tau \wedge T_0)} |\mu_0 \tau \wedge T_0 + \sigma W_{\tau \wedge T_0}|] \\ &\leq K_1 |x-y|.\end{aligned}$$

According to Corollary 3.7, the set $\{(\hat{V}_0^{(k)}(x), \hat{V}_1^{(k)}(x)), k \in \mathbb{N}\}$ is bounded for every $x > 0$. Therefore, the Arzela–Ascoli theorem gives that the increasing sequence $(\hat{V}_0^{(k)}, \hat{V}_1^{(k)})$ converges uniformly on every compact subset of $\mathbb{R}_+$ to some $(\hat{V}_0^{(\infty)}, \hat{V}_1^{(\infty)})$.

On the other hand, for a fixed $k$, $(\hat{V}_0^{(k)}, \hat{V}_1^{(k)})$ is the unique viscosity solution with linear growth to the system of variational inequalities

$$F_0^{(k)}(u_0, u_0', u_0'') = \min(\rho u_0 - \mathcal{L}_0 u_0, u_0 - \max(\hat{V}_0, \hat{V}_1^{(k-1)}(\cdot - g))) = 0,$$
$$F_1(u_1, u_1', u_1'') = \min(\rho u_1 - \mathcal{L}_1 u_1, u_1' - 1) = 0,$$

with initial condition $u_0(0) = 0$, $u_1(0) = \hat{V}_0^{(k)}((1-\lambda)g)$.

Since $\hat{V}_1^{(k-1)}$ converges uniformly on every compact subset of $\mathbb{R}_+$, the Hamiltonian $F_0^{(k)}$ converges to $F_0$ on every compact subset of $\mathbb{R} \times \mathbb{R} \times \mathbb{R}$, with

$$F_0(u, u', u'') = \min(\rho u - \mathcal{L}_0 u, u - \max(\hat{V}_0, \hat{V}_1^{\infty}(\cdot - g))) = 0.$$

According to standard stability results for viscosity solution see, for instance, Lemma 6.2, page 73, in Fleming and Soner [9], the couple $(\hat{V}_0^{(\infty)}, \hat{V}_1^{(\infty)})$ is a viscosity solution of the system of variational inequalities

(5.6) $\quad \min(\rho \hat{V}_0^{\infty} - \mathcal{L}_0 \hat{V}_0^{\infty}, \hat{V}_0^{\infty} - \max(\hat{V}_0, \hat{V}_1^{\infty}(\cdot - g))) = 0,$

(5.7) $\quad \min(\rho \hat{V}_1^{\infty} - \mathcal{L}_1 \hat{V}_1^{\infty}, (\hat{V}_1^{\infty})' - 1) = 0,$

with initial conditions $\hat{V}_1^{\infty}(0) = \hat{V}_0^{\infty}((1-\lambda)g)$ and $\hat{V}_0^{\infty}(0) = 0$. By uniqueness to the system (5.6)–(5.7), we conclude that $\hat{V}_0^{\infty} = \bar{v}_0 = v_0$ and $\hat{V}_1^{\infty} = v_1$. $\square$

We will close this section by describing the optimal strategy. According to Proposition 5.1, the value functions can be constructed recursively starting from $(\hat{V}_0, \hat{V}_1)$. Two cases have to be considered:

*Case* A: $\hat{V}_0^{(1)}((1-\lambda)g) = \hat{V}_0((1-\lambda)g)$. Then we have

$$\hat{V}_1^{(1)}(x) = \sup_{Z \in \mathcal{Z}} \mathbb{E}\left[\int_0^{T_1^-} e^{-\rho t} dZ_t + e^{-\rho T_1} \hat{V}_0^{(1)}((1-\lambda)g)\right]$$
$$= \sup_{Z \in \mathcal{Z}} \mathbb{E}\left[\int_0^{T_1^-} e^{-\rho t} dZ_t + e^{-\rho T_1} \hat{V}_0((1-\lambda)g)\right]$$
$$= \hat{V}_1(x).$$

Therefore, we deduce by a straightforward induction that the sequence $(\hat{V}_0^{(k)})_k$ is constant for $k \geq 1$ and the sequence $(\hat{V}_1^{(k)})_k$ is constant for $k \geq 0$. Therefore, we deduce from Proposition 5.1 that $v_0 = \hat{V}_0^{(1)}$ and $v_1 = \hat{V}_1$.

In regime 0, the optimal strategy consists in computing the optimal thresholds $a$ and $x_{01}^*$ associated to the optimal stopping problem $\hat{V}_0^{(1)}$. It is optimal to switch from regime 0 to regime 1 if the state process $R^0$ crosses the threshold $x_{01}^*$ while it is optimal to pay dividends and therefore abandon the growth opportunity forever if $R^0$ falls below the threshold $a$. At the level $a$, it is too costly to wait reaching the threshold $x_{01}^*$ even if the growth option is valuable. The shareholders prefer to receive today dividends rather than waiting for a more profitable payment in the future.

In regime 1, the optimal strategy consists in paying dividends above $\hat{x}_1$ and switching to regime 0 only when the manager is being forced by its cash constraints.

*Case* B: $\hat{V}_0^{(1)}((1-\lambda)g) > \hat{V}_0((1-\lambda)g)$. Let us introduce the sequence

$$\hat{\theta}_0^{(k)}(x) = \sup_{\tau \in \mathcal{T}} \mathbb{E}[e^{-\rho(\tau \wedge T_0)} \hat{\theta}_1^{(k-1)}(R_{\tau \wedge T_0}^{x,0} - g)], \qquad k \geq 1,$$

$$\hat{\theta}_1^{(k)}(x) = \sup_{Z \in \mathcal{Z}} \mathbb{E}\left[\int_0^{T_1^-} e^{-\rho t} dZ_t + e^{-\rho T_1} \hat{\theta}_0^{(k)}((1-\lambda)g)\right], \qquad k \geq 1,$$

starting from $\hat{\theta}_1^{(0)} = \hat{V}_1$ and $\hat{\theta}_0^{(0)} = \hat{V}_0$. Proceeding analogously as in the proof of Proposition 5.1, we can prove that the sequence $(\theta_0^{(k)}, \theta_1^{(k)})$ converges to $(\theta_0^{(\infty)}, \theta_1^{(\infty)})$ solution of the system of variational inequalities

$$\min(\rho \hat{\theta}_0^\infty - \mathcal{L}_0 \hat{\theta}_0^\infty, \hat{\theta}_0^\infty - \hat{\theta}_1^\infty(\cdot - g)) = 0,$$

$$\min(\rho \hat{\theta}_1^\infty - \mathcal{L}_1 \hat{\theta}_1^\infty, (\hat{\theta}_1^\infty)' - 1) = 0,$$

with initial conditions $\hat{\theta}_1^\infty(0) = \hat{\theta}_0^\infty((1-\lambda)g)$ and $\hat{\theta}_0^\infty(0) = 0$.

Note that the function $\hat{\theta}_0^\infty$ corresponds to the managerial decision to accumulate cash reserve at the expense of the shareholder's dividend payment in order to invest in the modern technology.

The key feature of our model in case B, which has to be viewed as the analogue of Proposition 3.5 in [5], can be summarized as follows:

If the net expected value evaluated at the threshold $\hat{x}_0$ dominates the firm value running under the old technology that is $\hat{\theta}_0^\infty(\hat{x}_0) > \hat{V}_0(\hat{x}_0)$, then the manager postpones dividend distribution in order to invest in the modern technology and, thus, $v_0 = \hat{\theta}_0^\infty$. Moreover, in regime 1, the manager always prefers to run under the modern technology until the cash process $X_t^1$ reaches zero, forcing the manager to return back in regime 0 with the value $\hat{\theta}_0^\infty((1-\lambda)g)$, that is, $v_1 = \hat{\theta}_1^\infty$.

If, on the contrary, $\hat{\theta}_0^\infty(\hat{x}_0) \leq \hat{V}_0(\hat{x}_0)$, then the manager optimally ignores the strategy $\hat{\theta}_0^\infty$. Several situations can occur. For small values of the cash process ($X_t^0 \leq a$), the manager optimally runs the firm under the old technology and pays out any surplus above $\hat{x}_0$ as dividends. For high values of the cash process ($X_t^0 \geq x_{01}^*$), the manager switches optimally to regime one. For intermediary values of the cash process ($a \leq X_t^0 \leq x_{01}^*$), there is an inaction region where the manager has not enough information to decide whether or not the investment is valuable.

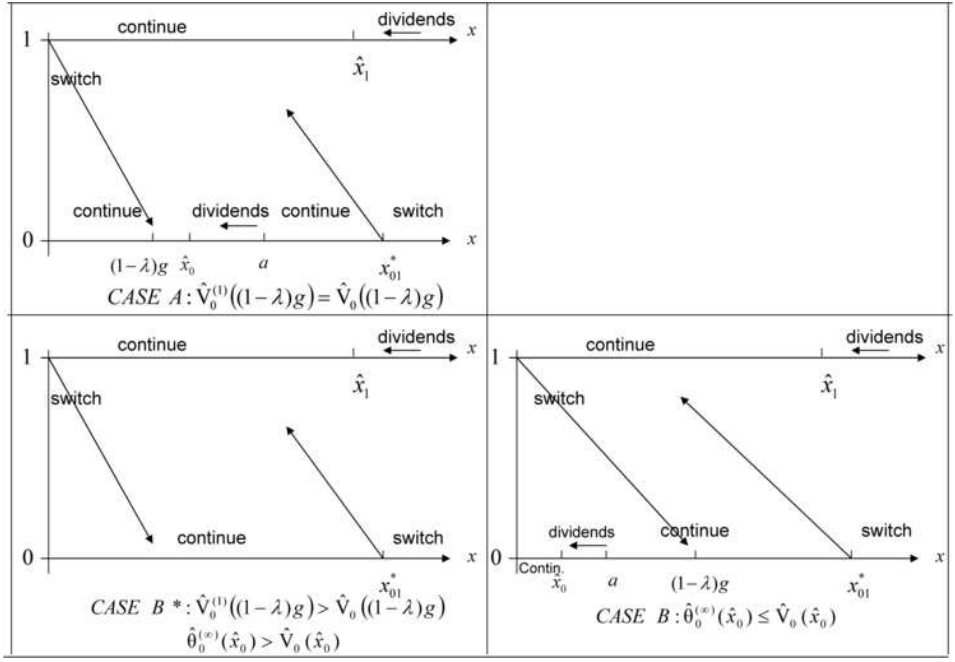

Fig. 1.

We summarize all the results in Synthetic Table 1 and Figure 1.

Synthetic Table 1

| $\frac{\mu_1}{\rho} \le \hat{V}_0((1-\lambda)g)$ | $\hat{V}_0((1-\lambda)g) \le \frac{\mu_1}{\rho} \le \frac{\mu_0}{\rho} + \hat{x}_1 + g - \hat{x}_0$ | $\frac{\mu_1}{\rho} > \max(\hat{V}_0((1-\lambda)g),$ $\frac{\mu_0}{\rho} + \hat{x}_1 - g - \hat{x}_0)$ |
|---|---|---|
| $v_0(x) = \hat{V}_0(x)$ | $v_0(x) = \hat{V}_0(x)$ | $v_0(x) = \hat{V}_0^\infty(x)$ |
| $v_1(x) = x + (1-\lambda)g - \hat{x}_0 + \frac{\mu_0}{\rho}$ | $v_1(x) = \hat{V}_1(x)$ | $v_1(x) = \hat{V}_1^\infty(x)$ |
| (diagram) | (diagram) | See Figure 1 |

5.3. *Computational aspects.* In the following lines, we briefly describe a way of computing the free boundary points that characterize the optimal stopping strategy.

*Case* A. In that case, the function $\hat{V}_1$ is given by formula (4.5) with $L = \hat{V}_0((1-\lambda)g)$ and the optimal threshold $\hat{x}_1$ is implicitly given by (4.6) that can be solved using a standard Newton method. To figure out the optimal thresholds $a$ and $x_{01}^*$ of regime 0, we have to solve the following system where $(a, x_{01}^*, B, C)$ are unknown:

$$\begin{cases} Be^{m_0^+ a} + Ce^{m_0^- a} = V_0(a), \\ Be^{m_0^+ x_{01}^*} + Ce^{m_0^- x_{01}^*} = \hat{V}_1(x_{01}^* - g), \\ m_0^+ Be^{m_0^+ a} + m_0^- Ce^{m_0^- a} = V_0'(a), \\ m_0^+ Be^{m_0^+ x_{01}^*} + m_0^- Ce^{m_0^- x_{01}^*} = \hat{V}_1'(x_{01}^* - g). \end{cases}$$

*Case* B. The computation of the optimal thresholds in that case is more involved and relies on the computation of the approximating thresholds associated to the sequence $(\hat{\theta}_0^{(k)}, \hat{\theta}_1^{(k)})$. At the time step $k$, the function $\hat{\theta}_1^{(k)}$ is given by formula (4.5) with $L = \hat{\theta}_0^{(k)}((1-\lambda)g)$. The optimal threshold $\hat{x}_1^{(k)}$ associated to the value function $\hat{\theta}_1^{(k)}$ is then given by (4.6). In regime 0, the optimal threshold $x_{01}^{(k)}$ associated to the value function $\hat{\theta}_0^{(k)}$ can be figured out using the smooth fit principle, that is,

$$\begin{cases} Ae^{m_0^+ x_{01}^{(k)}} = \hat{\theta}_1^{(k-1)}(x_{01}^{(k)}), \\ m_0^+ Ae^{m_0^+ x_{01}^{(k)}} = (\hat{\theta}_1^{(k-1)})'(x_{01}^{(k)}). \end{cases}$$

## APPENDIX A: PROOF OF THEOREM 3.1

We divide the proof into several steps.

PROOF OF THE CONTINUITY OF $v_1$ ON $(0, \infty)$. We prove that $v_1$ is continuous at any $y > 0$. We fix an arbitrary small $\varepsilon > 0$. Applying the dynamic programming principle (DP) to $v_1$, there exists a control $\alpha = (Z, (\tau_n)_{n \geq 1}) \in \mathcal{A}$ s.t.

$$v_1(y) - \frac{\varepsilon}{3} \leq \mathbb{E}\bigg[\int_0^{(\tau_1 \wedge T)^-} e^{-\rho t} dZ_t + e^{-\rho(\tau_1 \wedge T)}(v_1(X_T^{y,1}) 1_{T < \tau_1}$$

(A.1)
$$+ v_0(X_{\tau_1}^{y,1}) 1_{\tau_1 \leq T})\bigg],$$

$$= \mathbb{E}\bigg[\int_0^{(\tau_1 \wedge T)^-} e^{-\rho t} dZ_t + e^{-\rho(\tau_1 \wedge T)} v_0(X_{\tau_1}^{y,1}) 1_{\tau_1 \leq T}\bigg],$$

with $T = T^{y,1,\alpha}$ the bankruptcy time of the process $X^{y,1,\alpha}$, and since $v_1(X_T^{y,1}) = 0$ for $X_T^{y,1} < 0$.

For any $0 < x < y$, let $\theta = T^{y-x,1,\alpha}$ be the bankruptcy time of the process $X^{y-x,1,\alpha}$. We notice that $\theta \leq T$ and $X^{y-x,1,\alpha} = X^{y,1,\alpha} - x$ for all $0 < t < \theta \leq T$. Applying the dynamic programming principle (DP), we then have

$$v_1(y-x)$$
$$\geq \mathbb{E}\bigg[\int_0^{(\theta \wedge \tau_1)^-} e^{-\rho t} dZ_t$$
$$+ e^{-\rho(\theta \wedge \tau_1)}(v_1(X_\theta^{y-x,1})1_{\theta < \tau_1} + v_0(X_{\tau_1}^{y-x,1})1_{\tau_1 \leq \theta})\bigg]$$

(A.2)
$$\geq \mathbb{E}\bigg[\int_0^{(\theta \wedge \tau_1)^-} e^{-\rho t} dZ_t + e^{-\rho(\theta \wedge \tau_1)} v_0(X_{\tau_1}^{y-x,1})1_{\tau_1 \leq \theta}\bigg]$$

$$\geq \mathbb{E}\bigg[\int_0^{(\tau_1 \wedge T)^-} e^{-\rho t} dZ_t + e^{-\rho(\tau_1 \wedge T)} v_0(X_{\tau_1}^{y,1})1_{\tau_1 \leq T}\bigg]$$

$$- \mathbb{E}\bigg[\int_{\theta \wedge \tau_1}^{(T \wedge \tau_1)^-} e^{-\rho t} dZ_t\bigg]$$

$$+ \mathbb{E}[e^{-\rho(\theta \wedge \tau_1)} v_0(X_{\tau_1}^{y-x,1})1_{\tau_1 \leq \theta} - e^{-\rho(T \wedge \tau_1)} v_0(X_{\tau_1}^{y,1})1_{\tau_1 \leq T}].$$

Notice that $\theta \to T$ as $x$ goes to zero. Hence, by the continuity of $v_0$ and the dominated convergence theorem, one can find $0 < \delta_1 < y$ s.t. for $0 < x < \delta_1$:

(A.3) $\quad \mathbb{E}[e^{-\rho(\theta \wedge \tau_1)} v_0(X_{\tau_1}^{y-x,1})1_{\tau_1 \leq \theta} - e^{-\rho(T \wedge \tau_1)} v_0(X_{\tau_1}^{y,1})1_{\tau_1 \leq T}] \geq -\frac{\varepsilon}{3}.$

We also have

$$-\mathbb{E}\bigg[\int_{\theta \wedge \tau_1}^{(T \wedge \tau_1)^-} e^{-\rho t} dZ_t\bigg] \geq -\mathbb{E}[Z_{(T \wedge \tau_1)^-} - Z_{\theta \wedge \tau_1}].$$

From the dominated convergence theorem, one can find $0 < \delta_2 < y$ s.t. for $0 < x < \delta_2$:

(A.4) $$-\mathbb{E}\bigg[\int_{\theta \wedge \tau_1}^{(T \wedge \tau_1)^-} e^{-\rho t} dZ_t\bigg] \geq -\frac{\varepsilon}{3}.$$

Plugging inequalities (A.3) and (A.4) into (A.2), we obtain for $0 < x < \min\{\delta_1, \delta_2\}$

$$v_1(y-x) \geq \mathbb{E}\bigg[\int_0^{(\tau_1 \wedge T)^-} e^{-\rho t} dZ_t + e^{-\rho(\tau_1 \wedge T)} v_0(X_{\tau_1}^{y,1})1_{\tau_1 \leq T}\bigg] - \frac{2\varepsilon}{3}.$$

Using the inequality (A.1), and recalling that $v_1$ is nondecreasing, this implies

$$0 \leq v_1(y) - v_1(y-x) \leq \varepsilon,$$

which shows the left-continuity of $v_1$. By proceeding exactly in the same manner, we may obtain for a given $y > 0$ and any arbitrary $\varepsilon > 0$, the existence of $0 < \delta < y$ such that, for all $0 < x < \delta$,

$$0 \leq v_1(y+x) - v_1(y) \leq \varepsilon,$$

which shows the right-continuity of $v_1$. $\square$

PROOF OF SUPERSOLUTION PROPERTY. Fix $i \in \{0,1\}$. Consider any $\bar{x} \in (0, \infty)$ and $\varphi \in C^2(0, \infty)$ s.t. $\bar{x}$ is a minimum of $v_i - \varphi$ in a neighborhood $B_\varepsilon(\bar{x}) = (\bar{x} - \varepsilon, \bar{x} + \varepsilon)$ of $\bar{x}$, $\bar{x} > \varepsilon > 0$, and $v_i(\bar{x}) = \varphi(\bar{x})$. First, by considering the admissible control $\bar{\alpha} = (\overline{Z}, \bar{\tau}_n, n \geq 1)$ where we decide to take immediate switching control, that is, $\bar{\tau}_1 = 0$, while deciding not to distribute any dividend $\overline{Z} = 0$, we obtain

(A.5) $$v_i(\bar{x}) \geq v_{i-1}(\bar{x} - g_{i,1-i}).$$

On the other hand, let us consider the admissible control $\hat{\alpha} = (\hat{Z}, \hat{\tau}_n, n \geq 1)$ where we decide to never switch regime, while the dividend policy is defined by $\hat{Z}_t = \eta$ for $t \geq 0$, with $0 \leq \eta \leq \varepsilon$. Define the exit time $\tau_\varepsilon = \inf\{t \geq 0, X_t^{\bar{x},i} \notin \overline{B}_\varepsilon(\bar{x})\}$. We notice that $\tau_\varepsilon < T$. From the dynamic programming principle (DP), we have

(A.6)
$$\varphi(\bar{x}) = v(\bar{x}) \geq \mathbb{E}\left[\int_0^{\tau_\varepsilon \wedge h} e^{-\rho t} d\hat{Z}_t + e^{-\rho(\tau_\varepsilon \wedge h)} v_i(X_{\tau_\varepsilon \wedge h}^{\bar{x},i})\right]$$
$$\geq \mathbb{E}\left[\int_0^{\tau_\varepsilon \wedge h} e^{-\rho t} d\hat{Z}_t + e^{-\rho(\tau_\varepsilon \wedge h)} \varphi(X_{\tau_\varepsilon \wedge h}^{\bar{x},i})\right].$$

Applying Itô's formula to the process $e^{-\rho t} \varphi(X_t^{\bar{x},i})$ between 0 and $\tau_\varepsilon \wedge h$, and taking the expectation, we obtain

(A.7)
$$\mathbb{E}[e^{-\rho(\tau_\varepsilon \wedge h)} \varphi(X_{\tau_\varepsilon \wedge h}^{\bar{x},i})] = \varphi(\bar{x}) + \mathbb{E}\left[\int_0^{\tau_\varepsilon \wedge h} e^{-\rho t}(-\rho\varphi + \mathcal{L}_i\varphi)(X_t^{\bar{x},i})\, dt\right]$$
$$+ \mathbb{E}\left[\sum_{0 \leq t \leq \tau_\varepsilon \wedge h} e^{-\rho t}[\varphi(X_t^{\bar{x},i}) - \varphi(X_{t^-}^{\bar{x},i})]\right].$$

Combining relations (A.6) and (A.7), we have

(A.8)
$$\mathbb{E}\left[\int_0^{\tau_\varepsilon \wedge h} e^{-\rho t}(\rho\varphi - \mathcal{L}_i\varphi)(X_t^{\bar{x},i})\, dt\right] - \mathbb{E}\left[\int_0^{\tau_\varepsilon \wedge h} e^{-\rho t} d\hat{Z}_t\right]$$
$$- \mathbb{E}\left[\sum_{0 \leq t \leq \tau_\varepsilon \wedge h} e^{-\rho t}[\varphi(X_t^{\bar{x},i}) - \varphi(X_{t^-}^{\bar{x},i})]\right] \geq 0.$$

Take first $\eta = 0$. We then observe that $X$ is continuous on $[0, \tau_\varepsilon \wedge h]$ and only the first term of the relation (A.8) is nonzero. By dividing the above inequality by $h$ with $h \to 0$, we conclude that

(A.9) $$(\rho\varphi - \mathcal{L}_i\varphi)(\bar{x}) \geq 0.$$

Take now $\eta > 0$ in (A.8). We see that $\hat{Z}$ jumps only at $t = 0$ with size $\eta$, so that

(A.10) $$\mathbb{E}\left[\int_0^{\tau_\varepsilon \wedge h} e^{-\rho t}(\rho\varphi - \mathcal{L}_i\varphi)(X_t^{\bar{x},i})\,dt\right] - \eta - (\varphi(\bar{x} - \eta) - \varphi(\bar{x})) \geq 0.$$

By sending $h \to 0$, and then dividing by $\eta$ and letting $\eta \to 0$, we obtain

(A.11) $$\varphi'(\bar{x}) - 1 \geq 0.$$

This proves the required supersolution property

(A.12) $$\min[(\rho\varphi - \mathcal{L}_i\varphi)(\bar{x}), \varphi'(\bar{x}) - 1, v_i(\bar{x}) - v_{1-i}(\bar{x} - g_{i,1-i})] \geq 0. \qquad \square$$

PROOF OF THE SUBSOLUTION PROPERTY. We prove the subsolution property by contradiction. Suppose that the claim is not true. Then, there exists $\bar{x} > 0$ and a neighborhood $B_\varepsilon(\bar{x}) = (\bar{x} - \varepsilon, \bar{x} + \varepsilon)$ of $\bar{x}$, $\bar{x} > \varepsilon > 0$, a $C^2$ function $\varphi$ with $(\varphi - v_*)(\bar{x}) = 0$ and $\varphi \geq v_i$ on $\overline{B}_\varepsilon(\bar{x})$, and $\eta > 0$, s.t. for all $x \in \overline{B}_\varepsilon(\bar{x})$:

(A.13) $$\rho\varphi(x) - \mathcal{L}_i\varphi(x) > \eta,$$

(A.14) $$\varphi'(x) - 1 > \eta,$$

(A.15) $$v_i(x) - v_{i-1}(x - g_{i,1-i}) > \eta.$$

For any admissible control $\alpha = (Z, \tau_n, n \geq 1)$, consider the exit time $\tau_\varepsilon = \inf\{t \geq 0, X_t^{\bar{x},i} \notin \overline{B}_\varepsilon(\bar{x})\}$. We notice that $\tau_\varepsilon < T$. Applying Itô's formula to the process $e^{-\rho t}\varphi(X_t^{\bar{x},i})$ between $0$ and $(\tau_\varepsilon \wedge \tau_1)^-$, and by noting that before $(\tau_\varepsilon \wedge \tau_1)^-$, $X^{x,i}$ stays in regime $i$ and in the ball $\overline{B}_\varepsilon(\bar{x})$, we obtain

(A.16)
$$\begin{aligned}
\mathbb{E}[e^{-\rho(\tau_\varepsilon \wedge \tau_1)^-}&\varphi(X_{(\tau_\varepsilon \wedge \tau_1)^-}^{\bar{x},i})] \\
&= \varphi(\bar{x}) + \mathbb{E}\left[\int_0^{(\tau_\varepsilon \wedge \tau_1)^-} e^{-\rho t}(-\rho\varphi(X_t^{\bar{x},i}) + \mathcal{L}_i\varphi(X_t^{\bar{x},i}))\,dt\right] \\
&\quad - \mathbb{E}\left[\int_0^{(\tau_\varepsilon \wedge \tau_1)^-} e^{-\rho t}\varphi'(X_t^{\bar{x},i})\,dZ_t^c\right] \\
&\quad + \mathbb{E}\left[\sum_{0 \leq t < \tau_\varepsilon \wedge \tau_1} e^{-\rho t}[\varphi(X_t^{\bar{x},i}) - \varphi(X_{t^-}^{\bar{x},i})]\right].
\end{aligned}$$

From Taylor's formula and (A.14), and noting that $\Delta X_t^{\bar{x},i} = -\Delta Z_t$ for all $0 \leq t < \tau_\varepsilon \wedge \tau_1$, we have

$$\varphi(X_t^{\bar{x},i}) - \varphi(X_{t^-}^{\bar{x},i}) = \Delta X_t^{\bar{x},i} \varphi'(X_t^{\bar{x},i} + z\Delta X_t^{\bar{x},i}) \quad (A.17)$$
$$\leq -(1+\eta)\Delta Z_t.$$

Plugging the relations (A.13), (A.14) and (A.17) into (A.16), we obtain

$$v_i(\bar{x}) = \varphi(\bar{x}) \geq \mathbb{E}\left[\int_0^{(\tau_\varepsilon \wedge \tau_1)^-} e^{-\rho t}\, dZ_t + e^{-\rho(\tau_\varepsilon \wedge \tau_1)^-} \varphi(X_{(\tau_\varepsilon \wedge \tau_1)^-}^{\bar{x},i})\right]$$

$$+ \eta\left(\mathbb{E}\left[\int_0^{(\tau_\varepsilon \wedge \tau_1)^-} e^{-\rho t}\, dt\right] + \mathbb{E}\left[\int_0^{(\tau_\varepsilon \wedge \tau_1)^-} e^{-\rho t}\, dZ_t\right]\right)$$

$$(A.18) \qquad \geq \mathbb{E}\left[\int_0^{(\tau_\varepsilon \wedge \tau_1)^-} e^{-\rho t}\, dZ_t + e^{-\rho \tau_\varepsilon^-} \varphi(X_{\tau_\varepsilon^-}^{\bar{x},i}) 1_{\tau_\varepsilon < \tau_1}\right.$$

$$\left. + e^{-\rho \tau_1^-} \varphi(X_{\tau_1^-}^{\bar{x},i}) 1_{\tau_1 \leq \tau_\varepsilon}\right]$$

$$+ \eta\left(\mathbb{E}\left[\int_0^{(\tau_\varepsilon \wedge \tau_1)^-} e^{-\rho t}\, dt\right] + \mathbb{E}\left[\int_0^{(\tau_\varepsilon \wedge \tau_1)^-} e^{-\rho t}\, dZ_t\right]\right).$$

Notice that while $X_{\tau_\varepsilon^-}^{\bar{x},i} \in \overline{B}_\varepsilon(\bar{x})$, $X_{\tau_\varepsilon}^{\bar{x},i}$ is either on the boundary $\partial \overline{B}_\varepsilon(\bar{x})$ or out of $\overline{B}_\varepsilon(\bar{x})$. However, there is some random variable $\gamma$ valued in $[0,1]$ s.t.

$$X^{(\gamma)} = X_{\tau_\varepsilon^-}^{\bar{x},i} + \gamma \Delta X_{\tau_\varepsilon}^{\bar{x},i}$$
$$= X_{\tau_\varepsilon^-}^{\bar{x},i} - \gamma \Delta Z_{\tau_\varepsilon} \in \partial \overline{B}_\varepsilon(\bar{x}).$$

Then similarly as in (A.17), we have

$$(A.19) \qquad \varphi(X^{(\gamma)}) - \varphi(X_{\tau_\varepsilon^-}^{\bar{x},i}) \leq -\gamma(1+\eta)\Delta Z_{\tau_\varepsilon}.$$

Noting that $X^{(\gamma)} = X_{\tau_\varepsilon}^{\bar{x},i} + (1-\gamma)\Delta Z_{\tau_\varepsilon}$, we have

$$(A.20) \qquad v_i(X^{(\gamma)}) \geq v_i(X_{\tau_\varepsilon}^{\bar{x},i}) + (1-\gamma)\Delta Z_{\tau_\varepsilon}.$$

Recalling that $\varphi(X^{(\gamma)}) \geq v_i(X^{(\gamma)})$, inequalities (A.19) and (A.20) imply

$$\varphi(X_{\tau_\varepsilon^-}) \geq v_i(X_{\tau_\varepsilon}^{\bar{x},i}) + (1+\gamma\eta)\Delta Z_{\tau_\varepsilon}.$$

Plugging into (A.18) and using (A.15), we have

$$v_i(\bar{x}) \geq \mathbb{E}\left[\int_0^{(\tau_\varepsilon \wedge \tau_1)^-} e^{-\rho t}\, dZ_t + e^{-\rho \tau_\varepsilon} v_i(X_{\tau_\varepsilon}^{\bar{x},i}) 1_{\tau_\varepsilon < \tau_1}\right.$$

$$\left. + e^{-\rho \tau_1} v_{1-i}(X_{\tau_1}^{\bar{x},i}) 1_{\tau_1 \leq \tau_\varepsilon}\right]$$

$$\text{(A.21)} \qquad + \eta \mathbb{E}\left[ \int_0^{\tau_\varepsilon \wedge \tau_1} e^{-\rho t}\, dt + \int_0^{(\tau_\varepsilon \wedge \tau_1)^-} e^{-\rho t}\, dZ_t + e^{-\rho \tau_1} \mathbf{1}_{\tau_1 \leq \tau_\varepsilon} \right.$$

$$\left. + \gamma e^{-\rho \tau_\varepsilon \wedge \tau_1} \Delta Z_{\tau_\varepsilon} \mathbf{1}_{\tau_\varepsilon < \tau_1} \right]$$

$$+ \mathbb{E}[e^{-\rho \tau_\varepsilon} \Delta Z_{\tau_\varepsilon} \mathbf{1}_{\tau_\varepsilon < \tau_1}].$$

We now claim that there exists a constant $c_0 > 0$ such that, for any admissible control

$$\text{(A.22)} \qquad \mathbb{E}\left[ \int_0^{\tau_\varepsilon \wedge \tau_1} e^{-\rho t}\, dt + \int_0^{(\tau_\varepsilon \wedge \tau_1)^-} e^{-\rho t}\, dZ_t \right.$$

$$\left. + e^{-\rho \tau_1} \mathbf{1}_{\tau_1 \leq \tau_\varepsilon} + \gamma e^{-\rho \tau_\varepsilon \wedge \tau_1} \Delta Z_{\tau_\varepsilon} \mathbf{1}_{\tau_\varepsilon < \tau_1} \right] \geq c_0.$$

The $C^2$ function $\psi(x) = c_0[1 - \frac{(x-\bar{x})^2}{\varepsilon^2}]$, with

$$0 < c_0 \leq \min\left\{ \left(\rho + \frac{2}{\varepsilon}\mu_i + \frac{1}{\varepsilon^2}\sigma^2\right)^{-1}, \frac{\varepsilon}{2} \right\},$$

satisfies

$$\text{(A.23)} \quad \begin{cases} \min\{-\rho\psi + \mathcal{L}_i\psi + 1, 1 - \psi', -\psi + 1\} \geq 0, & \text{on } \overline{B}_\varepsilon(\bar{x}), \\ \psi = 0, & \text{on } \partial \overline{B}_\varepsilon(\bar{x}). \end{cases}$$

Applying Itô's formula, we then obtain

$$\text{(A.24)} \qquad 0 < c_0 = \psi(\bar{x}) \leq \mathbb{E}[e^{-\rho(\tau_\varepsilon \wedge \tau_1)} \psi(X^{\bar{x},i}_{(\tau_\varepsilon \wedge \tau_1)^-})]$$

$$+ \mathbb{E}\left[ \int_0^{\tau_\varepsilon \wedge \tau_1} e^{-\rho t} dt \right] + \mathbb{E}\left[ \int_0^{(\tau_\varepsilon \wedge \tau_1)^-} e^{-\rho t} dZ_t \right].$$

Noting that $\psi'(x) \leq 1$, we have

$$\psi(X^{\bar{x},i}_{\tau_\varepsilon^-}) - \psi(X^{(\gamma)}) \leq (X^{\bar{x},i}_{\tau_\varepsilon^-} - X^{(\gamma)}) = \gamma \Delta Z_{\tau_\varepsilon}.$$

Plugging into (A.24), we obtain

$$\text{(A.25)} \qquad 0 < c_0 \leq \mathbb{E}[e^{-\rho \tau_1} \psi(X^{\bar{x},i}_{\tau_1^-}) \mathbf{1}_{\tau_1 \leq \tau_\varepsilon}] + \mathbb{E}\left[ \int_0^{\tau_\varepsilon \wedge \tau_1} e^{-\rho t}\, dt \right]$$

$$+ \mathbb{E}\left[ \int_0^{(\tau_\varepsilon \wedge \tau_1)^-} e^{-\rho t}\, dZ_t \right] + \mathbb{E}[\gamma e^{-\rho \tau_\varepsilon} \Delta Z_{\tau_\varepsilon} \mathbf{1}_{\tau_\varepsilon < \tau_1}].$$

Since $\psi(x) \leq 1$ for all $x \in B_\varepsilon(\bar{x})$, this proves the claim (A.22).

Finally, by taking the supremum over all admissible control $\alpha$, and using the dynamic programming principle (DP), (A.21) implies $v_i(\bar{x}) \geq v_i(\bar{x}) + \eta c_0$, which is a contradiction. Thus, we obtain the required viscosity subsolution property

$$\text{(A.26)} \quad \min[(\rho\varphi - \mathcal{L}_i\varphi)(\bar{x}), \varphi'(\bar{x}) - 1, v_i(\bar{x}) - v_{i-1}(\bar{x} - g_{i,i-1})] \leq 0.$$

□

PROOF OF THE UNIQUENESS PROPERTY. Suppose $u_i$, $i = 0, 1$, are continuous viscosity subsolutions to the system of variational inequalities on $(0, \infty)$, and $w_i$, $i = 0, 1$, continuous viscosity supersolutions to the system of variational inequalities on $(0, \infty)$, satisfying the boundary conditions $u_i(0^+) \leq w_i(0^+)$, $i = 0, 1$, and the linear growth condition

$$|u_i(x)| + |w_i(x)| \leq C_1 + C_2 x \qquad \forall x \in (0, \infty), i = 1, 2, \tag{A.27}$$

for some positive constants $C_1$ and $C_2$. We want to prove that

$$u_i \leq w_i \qquad \text{on } (0, \infty), i = 0, 1.$$

*Step* 1. We first construct strict supersolutions to the system with suitable perturbations of $w_i$, $i = 0, 1$. We set

$$h_i(x) = A_i + B_i x + C x^2, \qquad x > 0,$$

where

$$A_0 = \frac{\mu_1 B_1 + C\sigma^2 + 1}{\rho} + \frac{C}{4}\left(\frac{B_1}{C} - 2\frac{\mu_1}{\rho}\right)^2$$

$$+ \frac{C}{4}\left(\frac{B_0}{C} - 2\frac{\mu_0}{\rho}\right)^2 + w_0(0^+) + w_1(0^+),$$

$$A_1 = A_0 + \frac{3}{2}g + \frac{g}{\lambda},$$

$$B_0 = 3, \qquad B_1 = 2 + \frac{2}{\lambda},$$

$$C = \frac{1}{\lambda g}.$$

We then define, for all $\gamma \in (0, 1)$, the continuous functions on $(0, \infty)$ by

$$w_i^\gamma = (1 - \gamma) w_i + \gamma h_i, \qquad i = 0, 1.$$

We then see that, for all $\gamma \in (0, 1), i = 0, 1$,

$$\begin{aligned}
w_i^\gamma(x) &- w_{1-i}^\gamma(x - g_{i,1-i}) \\
&= (1 - \gamma)[w_i(x) - w_{1-i}(x - g_{i,1-i})] \\
&\quad + \gamma[h_i(x) - h_{1-i}(x - g_{i,1-i})], \\
&\geq \gamma[(2C g_{i,1-i} + B_i - B_{1-i})x \\
&\qquad + A_i - A_{1-i} - C g_{i,1-i}^2 + B_{1-i} g_{i,1-i}], \\
&\geq \gamma \frac{g}{2}, \qquad i = 0, 1.
\end{aligned} \tag{A.28}$$

Furthermore, we also easily obtain

(A.29) $$h'_i(x) - 1 = B_i + 2Cx - 1 \geq 1.$$

A straight calculation will also provide us with the last required inequality, that is,

(A.30) $$\rho h_i(x) - \mathcal{L}_i h_i(x) \geq 1.$$

Combining (A.28), (A.29) and (A.30), this shows that $w_i^\gamma$ is a strict supersolution of the system: for $i = 0, 1$, we have on $(0, \infty)$

(A.31) $$\min[\rho w_i^\gamma(x) - \mathcal{L}_i w_i^\gamma(x), w_i^{\gamma'}(x) - 1, w_i^\gamma(x) - w_{i-1}^\gamma(x - g_{i,1-i})]$$
$$\geq \gamma \min\left\{1, \frac{g}{2}\right\} = \delta.$$

*Step* 2. In order to prove the comparison principle, it suffices to show that, for all $\gamma \in (0, 1)$,

$$\max_{i \in \{0,1\}} \sup_{(0,+\infty)} (u_i - w_i^\gamma) \leq 0,$$

since the required result is obtained by letting $\gamma$ to 0. We argue by contradiction and suppose that there exist some $\gamma \in (0, 1)$ and $i \in \{0, 1\}$, s.t.

(A.32) $$\theta := \max_{j \in \{0,1\}} \sup_{(0,+\infty)} (u_j - w_j^\gamma) = \sup_{(0,+\infty)} (u_i - w_i^\gamma) > 0.$$

Notice that $u_i(x) - w_i^\gamma(x)$ goes to $-\infty$ when $x$ goes to infinity. We also have $\lim_{x \to 0^+} u_i(x) - \lim_{x \to 0^+} w_i^\gamma(x) \leq \gamma(\lim_{x \to 0^+} w_i(x) - A_i) \leq 0$. Hence, by continuity of the functions $u_i$ and $w_i^\gamma$, there exists $x_0 \in (0, \infty)$ s.t.

$$\theta = u_i(x_0) - w_i^\gamma(x_0).$$

For any $\varepsilon > 0$, we consider the functions

$$\Phi_\varepsilon(x, y) = u_i(x) - w_i^\gamma(y) - \phi_\varepsilon(x, y),$$

$$\phi_\varepsilon(x, y) = \frac{1}{4}|x - x_0|^4 + \frac{1}{2\varepsilon}|x - y|^2,$$

for all $x, y \in (0, \infty)$. By standard arguments in the comparison principle, the function $\Phi_\varepsilon$ attains a maximum in $(x_\varepsilon, y_\varepsilon) \in (0, \infty)^2$, which converges (up to a subsequence) to $(x_0, x_0)$ when $\varepsilon$ goes to zero. Moreover,

(A.33) $$\lim_{\varepsilon \to 0} \frac{|x_\varepsilon - y_\varepsilon|^2}{\varepsilon} = 0.$$

Applying Theorem 3.2 in [4], we get the existence of $M_\varepsilon, N_\varepsilon \in \mathbb{R}$ such that

$$(p_\varepsilon, M_\varepsilon) \in J^{2,+} u_i(x_\varepsilon),$$
$$(q_\varepsilon, N_\varepsilon) \in J^{2,-} w_i^\gamma(y_\varepsilon)$$

and

(A.34) $$\begin{pmatrix} M_\varepsilon & 0 \\ 0 & N_\varepsilon \end{pmatrix} \leq D^2\phi_\varepsilon(x_\varepsilon, y_\varepsilon) + \varepsilon(D^2\phi(x_\varepsilon, y_\varepsilon))^2,$$

where

$$p_\varepsilon = D_x\phi_\varepsilon(x_\varepsilon, y_\varepsilon) = \frac{1}{\varepsilon}(x_\varepsilon - y_\varepsilon) + (x_\varepsilon - x_0)^3,$$

$$q_\varepsilon = -D_y\phi_\varepsilon(x_\varepsilon, y_\varepsilon) = \frac{1}{\varepsilon}(x_\varepsilon - y_\varepsilon),$$

$$D^2\phi_\varepsilon(x_\varepsilon, y_\varepsilon) = \begin{pmatrix} 3(x_\varepsilon - x_0)^2 + \frac{1}{\varepsilon} & -\frac{1}{\varepsilon} \\ -\frac{1}{\varepsilon} & \frac{1}{\varepsilon} \end{pmatrix}.$$

By writing the viscosity subsolution property of $u_i$ and the viscosity supersolution property (A.31) of $w_i^\gamma$, we have the following inequalities:

(A.35) $$\min\left\{ \rho u_i(x_\varepsilon) - \left(\frac{1}{\varepsilon}(x_\varepsilon - y_\varepsilon) + (x_\varepsilon - x_0)^3\right)\mu_i - \frac{1}{2}\sigma^2 M_\varepsilon, \right.$$
$$\left. \left(\frac{1}{\varepsilon}(x_\varepsilon - y_\varepsilon) + (x_\varepsilon - x_0)^3\right) - 1, u_i((x_\varepsilon) - u_{1-i}(x_\varepsilon - g_{i,1-i}) \right\} \leq 0,$$

(A.36) $$\min\left\{ \rho w_i^\gamma(y_\varepsilon) - \frac{1}{\varepsilon}(x_\varepsilon - y_\varepsilon)\mu_i - \frac{1}{2}\sigma^2 N_\varepsilon, \frac{1}{\varepsilon}(x_\varepsilon - y_\varepsilon) - 1, \right.$$
$$\left. w_i^\gamma(y_\varepsilon) - w_{i-1}^\gamma(x_\varepsilon - g_{i,1-i}) \right\} \geq \delta.$$

We then distinguish the following three cases:

*Case 1*: $u_i(x_\varepsilon) - u_{1-i}(x_\varepsilon - g_{i,1-i}) \leq 0$ in (A.35).
From the continuity of $u_i$ and by sending $\varepsilon \to 0$, this implies

(A.37) $$u_i(x_0) \leq u_{1-i}(x_0 - g_{i,1-i}).$$

On the other hand, from (A.36), we also have

$$w_i^\gamma(y_\varepsilon) - w_{i-1}^\gamma(x_\varepsilon - g_{i,1-i}) \geq \delta,$$

which implies, by sending $\varepsilon \to 0$ and using the continuity of $w_i$,

(A.38) $$w_i^\gamma(x_0) \geq w_{i-1}^\gamma(x_0 - g_{i,1-i}) + \delta.$$

Combining (A.37) and (A.38), we obtain

$$\theta = u_i(x_0) - w_i^\gamma(x_0) \leq u_{1-i}(x_0 - g_{i,1-i}) - w_{i-1}^\gamma(x_0 - g_{i,1-i}) - \delta,$$
$$\leq \theta - \delta,$$

which is a contradiction.

*Case* 2: $(\frac{1}{\varepsilon}(x_\varepsilon - y_\varepsilon) + (x_\varepsilon - x_0)^3) - 1 \leq 0$ in (A.35). Notice that, by (A.36), we have

$$\frac{1}{\varepsilon}(x_\varepsilon - y_\varepsilon) - 1 \geq \delta,$$

which implies in this case

$$(x_\varepsilon - x_0)^3 \leq -\delta.$$

By sending $\varepsilon$ to zero, we obtain again a contradiction.

*Case* 3: $\rho u_i(x_\varepsilon) - (\frac{1}{\varepsilon}(x_\varepsilon - y_\varepsilon) + (x_\varepsilon - x_0)^3)\mu_i - \frac{1}{2}\sigma^2 M_\varepsilon \leq 0$ in (A.35). From (A.36), we have

$$\rho w_i^\gamma(y_\varepsilon) - \frac{1}{\varepsilon}(x_\varepsilon - y_\varepsilon)\mu_i - \frac{1}{2}\sigma^2 N_\varepsilon \geq \delta,$$

which implies in this case

(A.39) $\quad \rho(u_i(x_\varepsilon) - w_i^\gamma(y_\varepsilon)) - \mu_i(x_\varepsilon - x_0)^3 - \frac{1}{2}\sigma^2(M_\varepsilon - N_\varepsilon) \leq -\delta.$

From (A.34), we have

$$\tfrac{1}{2}\sigma^2(M_\varepsilon - N_\varepsilon) \leq \tfrac{3}{2}\sigma^2(x_\varepsilon - x_0)^2[1 + 3\varepsilon(x_\varepsilon - x_0)].$$

Plugging it into (A.39) yields

$$\rho(u_i(x_\varepsilon) - w_i^\gamma(y_\varepsilon)) \leq \mu_i(x_\varepsilon - x_0)^3 + \tfrac{3}{2}\sigma^2(x_\varepsilon - x_0)^2[1 + 3\varepsilon(x_\varepsilon - x_0)] - \delta.$$

By sending $\varepsilon$ to zero and using the continuity of $u_i$ and $w_i^\gamma$, we obtain the required contradiction: $\rho\theta \leq -\delta < 0$. This ends the proof. □

## APPENDIX B: PROOF OF PROPOSITION 3.3

**$C^1$ property.** We prove in three steps that, for a given $i \in 0, 1$, $v_i$ is a $C^1$ function on $(0, \infty)$. Notice first that since $v_i$ is a strictly nondecreasing continuous function on $(0, \infty)$, it admits a nonnegative left and right derivative $v_i'^-(x)$ and $v_i'^+(x)$ for all $x > 0$.

*Step* 1. We start by proving that $v_i'^-(x) \geq v_i'^+(x)$ for all $x \in (0, \infty)$. Suppose, on the contrary, that there exists some $x_0$ such that $v_i'^-(x_0) < v_i'^+(x_0)$. Take then some $q \in (v_i'^-(x), v_i'^+(x))$, and consider the function

$$\varphi(x) = v_i(x_0) + q(x - x_0) + \frac{1}{2\varepsilon}(x - x_0)^2,$$

with $\varepsilon > 0$. Then $x_0$ is a local minimum of $v_i - \varphi_i$, with $\varphi'(x_0) = q$ and $\varphi''(x_0) = \frac{1}{\varepsilon}$. Therefore, we get the required contradiction by writing the supersolution inequality

$$0 \leq \rho v_i(x_0) - \mu_i\varphi'(x_0) - \frac{\sigma^2}{2}\varphi''(x_0) = \rho v_i(x_0) - \mu_i q - \frac{\sigma^2}{2\varepsilon},$$

and choosing $\varepsilon$ small enough.

*Step* 2. We now prove that, for $i=0,1$, $v_i$ is $C^1$ on $(0,\infty)\backslash\mathcal{S}_i$.

Suppose there exists some $x_0 \notin \mathcal{S}_i$ s.t. $v_i'^{-}(x_0) > v_i'^{+}(x_0)$. We then fix some $q \in (v_i'^{+}(x_0), v_i'^{-}(x_0))$ and consider the function

$$\varphi(x) = v_i(x_0) + q(x-x_0) - \frac{1}{2\varepsilon}(x-x_0)^2,$$

with $\varepsilon > 0$. Then $x_0$ is a local maximum of $v_i - \varphi$, with $\varphi'(x_0) = q > 1$, $\varphi''(x_0) = -\frac{1}{\varepsilon}$. Since $x_0 \notin \mathcal{S}_i$, the subsolution inequality property implies

$$\rho v_i(x_0) - \mu_i q + \frac{\sigma^2}{2\varepsilon} \leq 0,$$

which leads to a contradiction, by choosing $\varepsilon$ sufficiently small. By combining the results from step 1 and step 2, we obtain that $v_i$ is $C^1$ on the open set $(0,\infty)\backslash\mathcal{S}_i$.

*Step* 3. We now prove that $v_i$ is $C^1$ on $(0,\infty)$.

From step 2, we have to prove the $C^1$ property of $v_i$ on $\mathcal{S}_i^*$. Fix then some $x_0 \in \mathcal{S}_i^*$ so that $v_i(x_0) = v_{1-i}(x_0 - g_{i,1-i})$. Hence, $x_0$ is a minimum of $v_i - v_{1-i}(\cdot - g_{i,1-i})$, and so

(B.1) $\quad v_i'^{-}(x_0) - v_{1-i}'^{-}(x_0 - g_{i,1-i}) \leq v_i'^{+}(x_0) - v_{1-i}'^{+}(x_0 - g_{i,1-i}).$

Now, from Lemma 4.3, $x_0 - g_{i,1-i}$ belongs to the open set $(0,\infty)\backslash\mathcal{S}_{1-i}$. From step 2, $v_{1-i}$ is $C^1$ on $(0,\infty)\backslash\mathcal{S}_{1-i}$, and so $v_{1-i}'^{+}(x_0 - g_{i,1-i}) = v_{1-i}'^{-}(x_0 - g_{i,1-i})$. From (B.1), we thus obtain

$$v_i'^{-}(x_0) \leq v_i'^{+}(x_0),$$

which is the required result, since the reverse inequality is already satisfied from step 1.

**$C^2$ property.** We now turn to the proof of the $C^2$ property of $v_i$ on the open set $\mathcal{C}_i \cup \text{int}(\mathcal{D}_i)$ of $(0,\infty)$. Since it is clear that $v_i$ is $C^2$ on $\text{int}(\mathcal{D}_i)$ (where $v_i' = 1$), we only have to prove that $v_i$ is $C^2$ on $\mathcal{C}_i$. By standard arguments, we check that $v_i$ is a viscosity solution to

(B.2) $\quad \rho v_i(x) - \mathcal{L}_i v_i(x) = 0, \qquad x \in \mathcal{C}_i.$

Indeed, let $\bar{x} \in \mathcal{C}_i$ and $\varphi$ a $C^2$ function on $\mathcal{C}_i$ s.t. $\bar{x}$ is a local maximum of $v_i - \varphi$, with $v_i(\bar{x}) = \varphi(\bar{x})$. Then, $\varphi'(\bar{x}) = v_i'(\bar{x}) > 1$. By definition of $\mathcal{C}_i$, we also have $v_i(\bar{x}) > v_{1-i}(x - g_{i,1-i})$ and so from the subsolution viscosity property (A.26) of $v_i$, we have

$$\rho\varphi(\bar{x}) - \mathcal{L}_i\varphi(\bar{x}) \leq 0.$$

The supersolution inequality for (B.2) is immediate from (A.12).

Now, for any arbitrary bounded interval $(x_1, x_2) \subset \mathcal{C}_i$, consider the Dirichlet boundary linear problem:

$$\rho w(x) - \mathcal{L}_i w(x) = 0 \qquad \text{on } (x_1, x_2) \tag{B.3}$$

$$w(x_1) = v_i(x_1) \qquad w(x_2) = v_i(x_2). \tag{B.4}$$

Classical results provide the existence and uniqueness of a smooth $C^2$ function $w$ solution on $(x_1, x_2)$ to (B.3)–(B.4). In particular, this smooth function $w$ is a viscosity solution to (B.2) on $(x_1, x_2)$. From standard uniqueness results for (B.3)–(B.4), we get $v_i = w$ on $(x_1, x_2)$. From the arbitrariness of $(x_1, x_2) \subset \mathcal{C}_i$, this proves that $v_i$ is smooth $C^2$ on $\mathcal{C}_i$.

## REFERENCES

[1] BOGUSLAVSKAYA, E. (2003). On optimization of dividend flow for a company in a presence of liquidation value. Working paper. Available at http://www.boguslavsky.net/fin/index.html.
[2] BREKKE, K. and OKSENDAL, B. (1994). Optimal switching in an economic activity under uncertainty. *SIAM J. Control Optim.* **32** 1021–1036. MR1280227
[3] CHOULLI, T., TAKSAR, M. and ZHOU, X. Y. (2003). A diffusion model for optimal dividend distribution for a company with constraints on risk control. *SIAM J. Control Optim.* **41** 1946–1979. MR1972542
[4] CRANDALL, M., ISHII, H. and LIONS, P. L. (1992). User's guide to viscosity solutions of second order partial differential equations. *Bull. Amer. Math. Soc.* **27** 1–67. MR1118699
[5] DÉCAMPS, J. P. and VILLENEUVE, S. (2007). Optimal dividend policy and growth option. *Finance and Stochastics* **11** 3–27. MR2284010
[6] DAVIS, M. H. A. and ZERVOS, M. (1994). A problem of singular stochastic control with discretionary stopping. *Ann. Appl. Probab.* **4** 226–240. MR1258182
[7] DIXIT, A. and PINDICK, R. (1994). *Investment under Uncertainty*. Princeton Univ. Press.
[8] DUCKWORTH, K. and ZERVOS, M. (2001). A model for investment decisions with switching costs. *Ann. Appl. Probab.* **11** 239–260. MR1825465
[9] FLEMING, W. and SONER, M. (1993). *Controlled Markov Processes and Viscosity Solutions*. Springer, Berlin. MR1199811
[10] GUO, X. and PHAM, H. (2005). Optimal partially reversible investment with entry decision and general production function. *Stochastic Process. Appl.* **115** 705–736. MR2132595
[11] JEANBLANC, M. and SHIRYAEV, A. (1995). Optimization of the flow of dividends. *Russian Math. Survey* **50** 257–277. MR1339263
[12] KARATZAS, I., OCONE, D., WANG, H. and ZERVOS, M. (2000). Finite-fuel singular control with discretionary stopping. *Stochastics Stochastics Rep.* **71** 1–50. MR1813505
[13] LY VATH, V. and PHAM, H. (2007). Explicit solution to an optimal switching problem in the two-regime case. *SIAM J. Control Optim.* **46** 395–426. MR2309034
[14] PHAM, H. (2007). On the smooth-fit property for one-dimensional optimal switching problem. *Séminaire de Probabilités XL. Lecture Notes in Math.* **1899** 187–202. Springer, Berlin.

[15] RADNER, R. and SHEPP, L. (1996). Risk vs. profit potential: A model of corporate strategy. *J. Economic Dynamics and Control* **20** 1373–1393.
[16] SHREVE, S., LEHOCZKY, J. P. and GAVER, D. (1984). Optimal consumption for general diffusions with absorbing and reflecting barriers. *SIAM J. Control Optim.* **22** 55–75. MR0728672
[17] VILLENEUVE, S. (2007). On the threshold strategies and smooth-fit principle for optimal stopping problems. *J. Appl. Probab.* **44** 181–198. MR2312995
[18] ZERVOS, M. (2003). A problem of sequential entry and exit decisions combined with discretionary stopping. *SIAM J. Control Optim.* **42** 397–421. MR1982276

V. LY VATH
LABORATOIRE ANALYSE ET PROBABILITÉS
UNIVERSITÉ D'EVRY
BD. F. MITTERAND
91025 EVRY CEDEX
FRANCE
E-MAIL: lyvath@ensiee.fr

H. PHAM
LABORATOIRE DE PROBABILITÉS
   ET MODÈLES ALÉATOIRES
UNIVERSITÉ PARIS 7
175 RUE DU CHEVALERET
75013 PARIS
AND
INSTITUT UNIVERSITAIRE DE FRANCE
FRANCE
E-MAIL: pham@math.jussieu.fr

S. VILLENEUVE
TOULOUSE SCHOOL OF ECONOMICS (GREMAQ-IDEI)
UNIVERSITÉ DE TOULOUSE 1
MANUFACTURE DES TABACS
21 ALLÉE DE BRIENNE
31000 TOULOUSE
FRANCE
E-MAIL: stephane.villeneuve@univ-tlse1.fr